\documentclass[letterpaper]{article} 

\usepackage[usenames,dvipsnames]{color}
\usepackage{bm}
\usepackage{graphicx,fancyhdr,lscape,dsfont,amssymb}
\usepackage[fleqn]{amsmath}
\usepackage{amsxtra,amsfonts,latexsym,amssymb,euscript}
\usepackage{amstext}
\usepackage{color}
\usepackage{longtable}
\usepackage{psfrag}
\usepackage{subfigure}
\usepackage{multirow}
\usepackage{stmaryrd}
\usepackage[titletoc]{appendix}
\usepackage{natbib}
\usepackage{graphicx}
\usepackage{mathrsfs}
\usepackage{ulem}
\usepackage{mathrsfs}
\usepackage{bbm}
\usepackage{upgreek}
\usepackage[boxed]{algorithm2e}

\allowdisplaybreaks

\newcommand{\tensor}[1]{  {\bm {#1}} } 

\begin{document}

\title{A phase field approach for damage propagation in periodic microstructured materials}

\author{Francesca Fantoni$^{1,*}$, Andrea Bacigalupo$^{2,}$\footnote{Corresponding authors: Tel:+39 0583 4326613, $\hspace{10cm}$ 
 E-mail addresses: francesca.fantoni@unibs.it; andrea.bacigalupo@imtlucca.it} , Marco Paggi$^{2}$, Jos\'e Reionoso$^3$
\\
\begin{small}
$^{1}$ DICATAM, Universit\`a degli Studi di Brescia, via Branze 43, 25123, Brescia, Italy
\end{small}
\\
\begin{small}
$^{2}$ IMT School for Advanced Studies Lucca,  Piazza S.Francesco 19, 55100 Lucca, Italy
 \end{small}
 \\
\begin{small}
$^{3}$ Group of Elasticity and Strength of Materials, School of Engineering, University of Seville,
\end{small}
\\
\begin{small}
 Camino de los Descubrimientos s/n,
41092, Seville, Spain
 \end{small}
}

\maketitle

\begin{abstract}
In the present work, the evolution of damage in periodic composite materials is investigated through  a novel finite element-based multiscale computational approach. The methodology is developed by means of the original  combination of  homogenization methods with the phase field approach of fracture. This last is applied at the macroscale level on the equivalent homogeneous continuum, whose constitutive properties are obtained in closed form via a two-scale asymptotic homogenization scheme. 
The formulation allows considering different assumptions on the evolution of damage at the microscale (e.g., damage in the matrix and not in the inclusion/fiber), as well as the role played by the microstructural topology. Numerical results show that the proposed formulation leads to an apparent tensile strength and a post-peak branch of unnotched and notched specimens dependent not only on the internal length scale of the phase field approach, as for homogeneous materials, but also on the inclusion volumetric content and its shape. Down-scaling relations allow the full reconstruction of the microscopic fields at any point of the macroscopic model, as a simple post-processing operation.
\end{abstract}
\section{Introduction}
\label{Sec::Inroduction}
The increasing demand  in achieving lightweight structures with superior performance in terms of damage tolerance and load-bearing capacities has motivated the incorporation of composite materials in many different applications in the last decades. As a result, at present, composite-like structures can be found in biomechanics, renewable energy systems, aerospace and automotive sectors, to quote some practical fields with strong societal impact. The irruption of new manufacturing techniques has encouraged this trend with renovated interest, allowing  the production  of tailor-made composite materials with unprecedent capabilities by means of controlling the microstructural arrangement   to attain   the specific practical needs.

In this context, it is well established that the versatility and applicability of composite materials strongly depend upon of their inherent failure mechanisms that take place at different scales.  These failure events can be  generally categorized as     meso-and-macro-scale mechanisms  \citep{maimi2007897,REINOSO201737} and micro-scale failure modes \citep{HERRAEZ2015196,arteiro2015}, with a cumbersome effective connection across the scales. From a mechanical point of view, the heterogenous character of composite materials induces their  inherent anisotropic response at the structural level,  being this behavior hardly predictable and thus representing one of the major drawbacks for their practical use. In order to overcome current limitations,  the potential routes for achievement  of superior responses in next generation of composite materials would require  a thorough understanding of the influence of the micro-scale arrangements on    macroscopic responses,  in terms of constituents properties   and their   spatial distribution within the material,   and the reliable prediction of failure events at such materials.

Investigation of composite materials having periodic or quasi-periodic arrangements resorts  micromechanical approaches, which  describe the  material behavior in detail, but in general result in labor-intensive analyses. Multiscale techniques, based on homogenization approaches, could conveniently overcome these drawbacks, allowing to gather  a synthetic but accurate description of the heterogeneous material behavior  both for static and dynamic problems. Homogenization methods take into account the effects of the microscopic phases on the overall constitutive properties of such composites, also in the presence of multi-field phenomena.  
Specifically,
homogenization methods allow replacing an heterogeneous material with a homogeneous equivalent one that can be modeled through either a Cauchy or a nonlocal continuum.
Homogenization techniques have been a matter of intensive research within the last decades.  In general sense, the local and/or nonlocal homogenization methods can be classified with respect to the underlying fundamental hypotheses as: the asymptotic techniques \citep{Bensoussan1978,Bakhvalov1984,GambinKroner1989,Hubert1992,Allaire1992,Meguid1994,Boutin1996,
AndrianovBolshakov2008,Panasenko2009,Tran2012,Bacigalupo2014,fantoni2017multi,fantoni2018design}, the variational-asymptotic techniques \citep{Willis1981,Smyshlyaev2000,Smyshlyaev2009,BacigalupoGambarotta2014b,BacigalupoGambarotta2014,
BacigalupoMoriniPiccolroaz2014,DelToro2019} and many identification approaches, involving the analytical \citep{Sevostianov2005,Bigoni2007,BaccaBigoniDalCorsoVeber2013a,BaccaBigoniDalCorsoVeber2013b,
BaccaDalCorsoVeberBigoni2013,BacigalupoGambarotta2013,Sevostianov2013,Rizzi2019,Rizzi2019b}, and the computational techniques \citep{Forest1998,Ostoja1999,KouznetsovaGeersBrekelmans2002,Forest2002,Feyel2003,KouznetsovaGeers2004,
Kaczmarczyk2008,Yuan2008,BacigalupoGambarotta2010,DeBellis2011,ForestTrinh2011,AdessiDeBellisSacco2013,
Trovalusci2015,Otero2018,Dirrenberger2019}.

Regarding the prediction of failure events in engineering materials and structures, the advent of new computational capabilities has promoted the generation of different numerical tools including diffusive crack methods \citep{BAZANT-JIRASEK,comi1999,peerlings20017723,dimitrijevic20111199}, strong discontinuity procedures \citep{MoesEtAl99,linder20071391,oliver20067093} and cohesive-like crack approaches \citep{CAMACHO,ORTIZPANDOLFI,PW12,turon2018506},  among many others, where most of them rely on the exploitation of finite element(FE)-based procedures. Recent variational formulations and crack tracking algorithms based on the analogy between linear elastic fracture mechanics and standard dissipative systems can be found in \citep{SalvadoriFantoniJMPS2016,Salvadori2019}. Derived from its versatility for the estimation of   failure mechanisms due to crack initiation and growth, the seminal variational approach of fracture developed by  \cite{FRANCFORTMARIGOJMPS1998}, being denominated as the phase field approach of fracture, endows a smeared crack idealization that permits overcoming most of the limitations  of alternative numerical methods. This methodology attains a regularized modeling of Griffith-like   fracture \citep{griffith1921}  and can be conceived as a nonlocal damage method, which in its original format is especially suitable for triggering fracture in brittle materials. In this concern,  \cite{BourdinFrancfortMarigo,Bourdin2000} comprehensively developed the corresponding numerical treatment and  extended the applicability of the phase field method for different applications (see also \cite{Pham2013147,Pahm2011618}), providing a robust   implementation. Alternatively,   \cite{miehe2010} came up with  a rigorous analysis  with regard to the thermodynamic considerations within a engineering point of view.    One of the main ingredients for the outstanding progresses on phase field methods stems from the fact that this variational approach does offer very appealing aspects and can  easily be implemented into  multi-field finite element frameworks.    In view of the strong potential of the phase field methods, recent developments   encompassed its application to cohesive-like fracture \citep{verhoosel2013}, coupled damage-plasticity   \citep{Ambati2015,MIEHE2015486,MIEHE20161},  shells \citep{mieheIJNME2014,Areias2016b,Reinoso2017Shell},   thermo-elastic \citep{MIEHE2015449} and hydrogen embrittlement \citep{Martinez2018} applications,  defining alternative degradation functions \citep{wu201772,Sagrado2017}, among many others. Owing to its modular formulation, the   phase-field approach to fracture has proven to be a   powerful tool for  fracture characterization of many different materials such as   arterial walls \citep{GULTEKIN201823}, and anisotropic media \citep{Teichtmeister201711,Bleyer2018213,QUINTANASCOROMINAS2019899}, to quote a few of them. Additional  contributions in the field have addressed the accuracy of phase field methods evaluating the   reliability  of  different implementation schemes with special attention on the local irreversibility constraint \citep{LINSE2017307}. Recently, the conjecture of crack nucleation in brittle materials discussed in \citep{tanne2018} has provided a further insight into the physical interpretation to the length scale in phase field formulations $\ell$, leading to a direct link  of this variable with respect the material strength $\sigma_c$, see \citep{yvo3}. 

With regard to heterogeneous media, being of special interest for composites,   new developments  of phase field methods    generally accounted for the combination of bulk fracture with cracking along the existing interfaces. Specifically, \cite{yvo1}   employed a level set function to describe interface displacement jump using an unique damage variable for both fracturing events, i.e. bulk and interface cracking. Conversely,  \cite{KHISAMITOV2018452}  proposed a split of the dissipative term in the variational formalism distinguishing between bulk and interface fracture mechanisms, whereas an alternative formulation was developed in \citep{HANSENDORR201925}. The  distinction  between bulk and interface energy dissipation was seminally exploited by the authors through the coupling of the phase field   and cohesive-like methods for triggering bulk and interface cracking, respectively, being denominated as the PF-CZM  \citep{paggi2017}. The robustness of this methodology has been assessed by its successful application  to poly-crystalline materials \citep{PAGGI2018123}, layered ceramics \citep{CAROLLO20182994} and micro-mechanics of  fiber-reinforced composites \citep{teresa} in comparison with the so-called coupled criterion \citep{manticgarcia2012}.

In this respect, though there exist different numerical methods that allow reliable characterization of complex fracture phenomena at   macro- and micro-scales, to the best authors' knowledge,  the  coupling between both scales of observations for the efficient material tailoring has not been addressed so far with the exploitation of the phase field approach of fracture. Consequently, the primary goal of the current study is to establish a novel numerical methodology which links a variational phase field approach of fracture at the macro-scale with a homogenization-based techniques in order to account for the micro-structural information and the spread of damage depending on the material microstructure geometry and properties. This framework represents an interesting development towards the efficient multi-scale modeling of damage events in reinforced composite through the combination of the phase field method of fracture and a careful microstructure treatment by means of asymptotic homogenization techniques.

 The manuscript is arranged  as follows. Section \ref{Sec::Synopsis} describes the developed methodology. The treatment of the corresponding material modeling at both scales is outlined in Section \ref{Sec::AsymHom}. The description of the damage evolution at the macro-scale is presented in Section \ref{Sec::phase_field}. The proposed methodology is examined in Section \ref{Sec::Benchmark} for different benchmark applications, whereas the main conclusions are given in Section \ref{Sec::Conclusions}.

\section{Synopsis of the proposed approach}
\label{Sec::Synopsis}
\begin{figure}[h!]
  \centering
  \includegraphics[width=16cm]{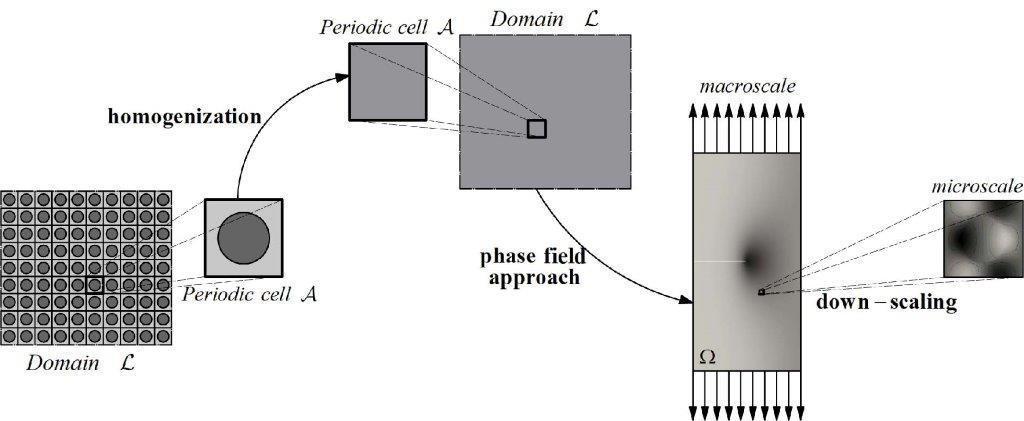}
  \caption{\it  Asymptotic homogenization allows describing the behavior of a microstructured periodic medium having domain $\mathcal{L}$ and periodic cell $\mathcal{A}$ as an equivalent homogeneous continuum, for different values of the phase field damage variable $\mathfrak{d}$. The evolution of the phase field is described by the phase field approach at the homogenized macroscale. Mechanical fields at the microscale can be eventually reconstructed through down-scaling relations.}
  \label{Fig::Synopsis}
\end{figure}

The present section highlights the key steps of the   methodology herein proposed for multi-scale simulation of damage within heterogeneous materials by originally combining asymptotic homogenization and the phase field (nonlocal) approach to damage, as sketched in Fig. \ref{Fig::Synopsis}.
In this respect, let consider a microstructured periodic medium, which is made by different phases and it is characterized by the macroscopic length scale $L$. At the microscale, each periodic cell is characterized by a microscopic length scale $\varepsilon$ which, for the validity of scale separation, has to be much smaller than $L$, i.e. $L\gg\varepsilon$.

Such a medium can be properly described trough a multiscale homogenization technique, which allows overcoming the computation cost of explicitly discretizing the heterogeneous microstructure using the finite element method. Thus, through the advocation of such a homogenization method, a concise but accurate description of material behavior and its mechanical damage can be gained. In this study, although the formulation could also be extended to other degradation mechanisms, we restrict our analysis to a situation where damage takes place only in the matrix, while the inclusions/fibers are assumed to behave obeying a linear elastic response with no degradation. This is a reasonable assumption in the majority of engineering applications concerning linear elastic stiff particle or fiber-reinforced composite materials. Hence, damage at the microscale is herein modeled as the  degradation of the elastic properties of the matrix.
 
A two-scale asymptotic homogenization procedure is therefore advocated in order to homogenize the microstructured composite medium for any damage level, whose phases are modeled as first-order elastic continua, to an equivalent first-order elastic continuum.

At the macroscale, specimens of any geometry and subjected to any loading conditions can be considered by employing the finite element method and the phase field approach to (nonlocal) damage propagation. Usually, the damaged constitutive tensor entering the stress-strain relation is degraded in the phase field approach according to a prescribed degradation function $g(\mathfrak{d})$ which does not take into account the damage mechanism occurring at the microscale. Following a standard procedure with the context of Damage Mechanics, depending on the phase field damage variable $\mathfrak{d}$, the elastic parameters of the matrix will be degraded by the degradation function $g(\mathfrak{d})=\left(1-\mathfrak{d}\right)^2$, while that of the particle/fiber will be left undamaged. This will establish the coupling with asymptotic homogenization, which is used to compute the degraded constitutive tensor based on the above assumption on the damage mechanism. As a result, the overall homogenized damaged constitutive tensor will not be simply equal to the undamaged one re-scaled by $g(\mathfrak{d})$ as in the standard phase field approach, but it will depend on the damage mechanisms at the microscale and eventually on the microscale topology (inclusion shape and volumetric content) and linear elastic parameters of the material constituents.

In this computational framework, from the practical standpoint, the asymptotic homogenization scheme is be called at each Gauss point of the finite element discretization at the macroscale, in order to compute the homogenized constitutive tensor of the damaged material, based on the damage-like phase field variable  $\mathfrak{d}$.
This resembles the FE$^2$ method, where two nested finite element models are solved, one at the macroscale, and another at the microscale. To reduce the computation cost associated with the present formulation, we herein explore a look-up table concept. Basically, the current asymptotic homogenization is applied off-line to a series of unit periodic cell problems for a series of values of $\mathfrak{d}$ ranging from zero to unity, i.e intact and fully degraded states, respectively. This information is subsequently employed to degrade the elastic properties of the matrix. The obtained homogenized constitutive tensor components, for different values of $\mathfrak{d}$, are then interpolated to determine their closed-form (analytical) dependency upon $\mathfrak{d}$, that to be used by the phase field approach for the computation of the damaged elastic strain energy density of the damaged material. Similarly, the first derivative of the components of the homogenized damaged constitutive tensor with respect to $\mathfrak{d}$ are also analytically determined, to be inserted into the computation of the tangent constitutive operator  at the macroscale, which is  required for fully consistent implicit iterative-incremental solution schemes.

Finally, as a post-process of the simulations, the homogenized stress and strain fields can be visualized. Moreover, by exploiting the down-scaling relations established by the asymptotic homogenization approach, the microscale fields can also be reconstructed and inspected, for any material point. Algorithm \ref{algorith} outlines the complete numerical-analytical approach herein developed for a given pseudo time increment throughout the simulation.

\begin{algorithm}[ht!]
\SetVlineSkip{8pt}
Initialization of the   displacement $\mathbf{u}(\mathbf{x})$ and   the fracture phase field variable $\mathfrak{d}(\mathbf{x})$  at the time $t_{n}$

 Within the pseudo time interval $\left[t_{n},t_{n+1}^{(k)} \right]$ at the iteration $k$, update the prescribed loading and the primary fields
$\{ \mathbf{u}_{n+1}^{(k)}(\mathbf{x}), \mathfrak{d}_{n+1}^{(k)}(\mathbf{x})  \}$  

 Loop over the integration points at the macroscale
\BlankLine
\begin{enumerate} 
\item Interpolate the macro-scale phase field variable $\overline{\mathfrak{d}}_{n+1}^{(k)}(\mathbf{x})$
   \item Call the asymptotic homogenization procedure for microscale characterization
	\item  Compute the homogenized constitutive tensor $C_{ijkl} (\mathfrak{d})$ (for the  unit cell) based on the look-up \\table concept  
\item Determine the first  derivative of $C_{ijkl} (\mathfrak{d})$ with respect to  $\mathfrak{d}$. 
\item Return $C_{ijkl} (\mathfrak{d})$ and $\partial_{\mathfrak{d}} C_{ijkl} (\mathfrak{d})$ to the macroscale
		\end{enumerate}
 Construct the element residual and stiffness operators.\\
 Solve for the new displacement increment 
	 \caption{Algorithm flowchart for the phase field-asymptotic homogenization procedure.}
	 \label{algorith}
\end{algorithm}

In summary, the main goal of the proposed formulation is to study the evolution of damage inside the bulk depending on the microstructural features of the material itself. In particular, by stating that damage takes place only in the matrix, the proposed methodology allows investigating the role of the material microstructure (inclusion shape and volumetric content) onto the macroscopic response, both in terms of load-displacement relation and in terms of damage propagation. Moreover, it is worth to state that the proposed methodology precludes the use of computational demanding FE$^2$ schemes, and whose use can be widely exploited in different applications.

A range of values for the characteristic phase field length scale $\ell$ are also examined, under the assumption that diffuse damage is spread enough in the space to avoid strain localization and therefore guarantee the applicability of homogenization theory.

As shown in the results obtained from the numerical simulations, the proposed coupled model is able to predict an apparent strength of the system that is not only affected by $\ell$ (which is usually linked to the apparent material strength), as in the standard phase field approaches for homogenous materials that do not account for different degradation scenarios at the microscale, but also on the microstructural topology and the inclusion volumetric content.

\section{Periodic elastic material modeled at two scales}
\label{Sec::AsymHom}
\begin{figure}[h!]
  \centering
  \includegraphics[width=10cm]{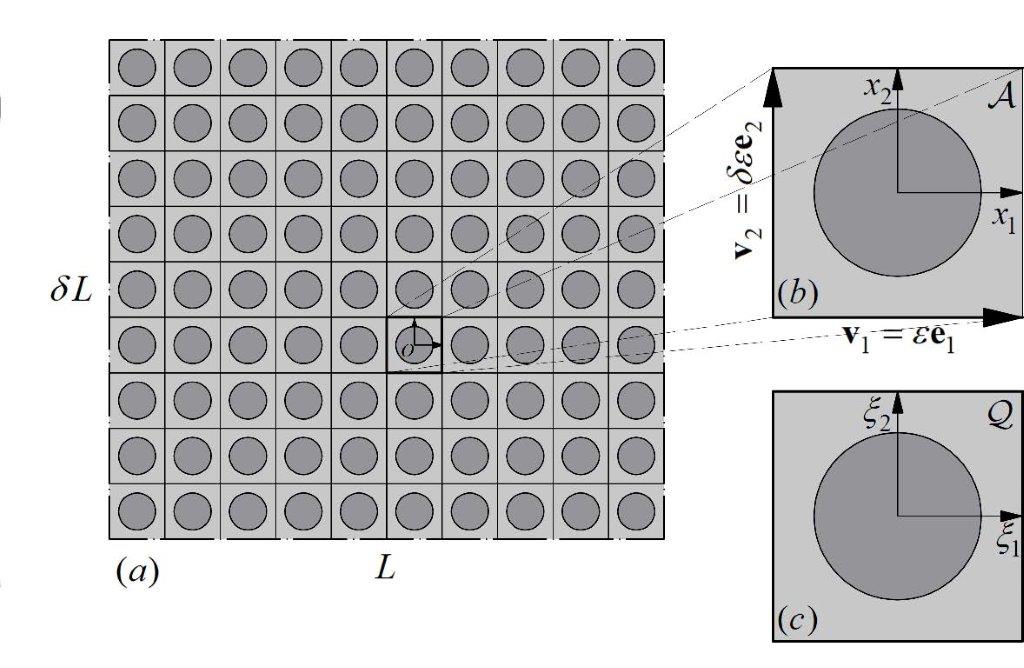}
  \caption{\it  (a) Heterogeneous microstructured medium having structural characteristic size $L$.
  (b) Periodic cell $\mathcal{A}$ with micro characteristic size $\varepsilon$ and periodicity vectors $\mathbf{v}_1$ and $\mathbf{v}_2$. (c) Periodic unit cell $\mathcal{Q}$. }
  \label{Fig::Homogenization}
\end{figure}

To establish a general scheme, one considers a linear elastic material characterized by a periodic  microstructure, whose phases can be described as first-order continua (Fig. \ref{Fig::Homogenization}).
In a two-dimensional perspective, vector $\mathbf{x}=x_1 \mathbf{e}_1 + x_2 \mathbf{e}_2$ identifies the position of each material point in the orthogonal coordinate system with origin at point $O$ and base $\mathbf{e}_1,\mathbf{e}_2$ as shown in Fig. \ref{Fig::Homogenization}-(b).

Advocating the periodicity of the medium microstructure, a periodic cell $\mathcal{A}=[0,\varepsilon]\times[0,\delta\varepsilon]$ is identified by two orthogonal periodicity vectors $\mathbf{v}_1=\varepsilon\mathbf{e}_1$ and $\mathbf{v}_2=\delta\varepsilon\,\mathbf{e}_2$, where $\varepsilon$ characterizes the size of cell $\mathcal{A}$ as depicted in Fig. \ref{Fig::Homogenization}.
As is usually done in the context of asymptotic homogenization, the unit cell $\mathcal{Q}=[0,1]\times[0,\delta]$ is through obtained rescaling the periodic cell $\mathcal{A}$ by the characteristic length $\varepsilon$.
The microscopic displacement $\mathbf{u}(\mathbf{x})=u_i(\mathbf{x})\mathbf{e}_i$ is   induced by body forces $\mathbf{b}(\mathbf{x})=b_i(\mathbf{x})\mathbf{e}_i$ whose periodicity is much greater than the microstructural one for the validity of the scale separation condition. Body forces in fact, are here assumed to be $\mathcal{L}$-periodic, and to have vanishing mean values on $\mathcal{L}$, being $\mathcal{L}=[0,L]\times[0,\delta L]$ a true representative portion of the whole body. Consistently with what above said, the structural length $L$ has to be much greater than the microstructural one (${L}\gg\varepsilon$). Moreover, diffuse damage without strain localization phenomena is also required to allow the rigorous applicability of homogenization theory.

Under the hypothesis of  quasi-static loading processes, the partial differential equation governing the elastic problem is given componentwise as
\begin{equation}
\label{eq:field_equation}
\frac{D}{D x_j}
\left(
C_{ijkl}^{(m,\varepsilon)}
\frac{D u_k}{D x_l}
\right)
+ b_i=0,
\end{equation}
being $C_{ijkl}^{(m,\varepsilon)}$ the components of the fourth order micro elasticity tensor, where the superscript $m$ refers to the microscale and $\varepsilon$ to the characteristic size of $\mathcal{A}$. In Eq. \eqref{eq:field_equation}, the symbol $D/D x_j$ denotes the generalized derivative with respect to the variable $x_j$.
The micro-constitutive elastic tensor is $\mathcal{Q}$-periodic, meaning that
\begin{equation}
\label{eq:micro_const_tensors}
C_{ijkl}^{(m,\varepsilon)}
\left(
\mathbf{x}
+
\mathbf{v}_{\alpha}\right)=
C_{ijkl}^{(m,\varepsilon)}
\left(
\mathbf{x}
\right),
\hspace{0.2 cm} \alpha=1,2, \hspace{0.2 cm}\forall \mathbf{x}\in\mathcal{A}
\end{equation}
and its components depend only on the fast  variable $\tensor{\xi}=\mathbf{x}/\varepsilon$:
\begin{equation}
\label{eq:Q_periodicity}
C_{ijkl}^{(m,\varepsilon)}
\left(
\mathbf{x}
\right)
=
C_{ijkl}^{m}
\left(
\tensor{\xi}=\frac{\mathbf{x}}{\varepsilon}
\right), \hspace{0.2 cm}\forall \tensor{\xi}\in \mathcal{Q}.
\end{equation}
Two variables in fact can be identified for the separation of the macro and the micro scales, namely the macroscopic (slow) one $\mathbf{x}\in\mathcal{A}$ and the microscopic (fast) one $\tensor{\xi}\in\mathcal{Q}$.
In this regard, the micro displacement results to be a function of both the slow and the fast variables, i.e. $u_i=u_i(\mathbf{x},\tensor{\xi}=\frac{\mathbf{x}}{\varepsilon})$.
In such a context, homogenization techniques reveal to be efficient in accurately describing the overall behavior of a microstructured composite material, replacing the periodic continuum with an homogeneous equivalent one requiring a much lower computation cost for its analysis and simulation.
Analytical and numerical solutions of Eq. \eqref{eq:field_equation} could be particularly complex to obtain because of the rapidly oscillating $\mathcal{Q}$-periodic coefficients. In what follows, the formulation of an equivalent first-order elastic continuum is described, together with the exact closed-form expression of the overall constitutive elastic tensor components.

In the obtained equivalent homogenized medium, the macro displacement field in each material point $\mathbf{x}$ is denoted with  $\mathbf{U}(\mathbf{x})=U_i(\mathbf{x})\mathbf{e}_i$.
It is worth mentioning that variability of source terms acting on the medium has to be much greater than the microstructural length scale $\varepsilon$, in order to preserve the scales separability principle. If volume forces are $\mathcal{L}$-periodic, then macroscopic displacement $\mathbf{U}(\mathbf{x})$ will be $\mathcal{L}$-periodic too, otherwise suitable boundary conditions should  be taken into account in order to determine the macroscopic field \citep{fantoni2017multi,fantoni2018design}.

\subsection{Asymptotic expansion of the microscopic field equation}
\label{SubSec:AsymExp}
%
%
%
Splitting the contributions of the slow and fast variables in accordance with the scale separation condition, the micro displacement field  can be asymptotically expanded in terms of the micro characteristic length scale $\varepsilon$ \citep{Bakhvalov1984} as
\begin{equation}
\label{eq:asymp_expansion}
u_k
\left(
\mathbf{x},\frac{\mathbf{x}}{\varepsilon}
\right)
=
\sum_{l=0}^{+\infty}\varepsilon^l \,u_k^{(l)}
\left(
\mathbf{x},
\frac{\mathbf{x}}{\varepsilon}
\right)=u_k^{(0)}
\left(
\mathbf{x},
\frac{\mathbf{x}}{\varepsilon}
\right)
+\varepsilon\,
u_k^{(1)}
\left(
\mathbf{x},
\frac{\mathbf{x}}{\varepsilon}
\right)
+\varepsilon^2\,
u_k^{(2)}
\left(
\mathbf{x},
\frac{\mathbf{x}}{\varepsilon}
\right)
+
\mathcal{O}(\varepsilon^3)
\end{equation}
Taking into account the differentiation rule for a generic function $f(\mathbf{x},\tensor{\xi}=\frac{\mathbf{x}}{\varepsilon})$:
\begin{equation}
\label{eq:derivation_rule}
\frac{D}{D x_j}f\left(
\mathbf{x},\tensor{\xi}=
\frac{\mathbf{x}}{\varepsilon}
\right)
=
\left.
\left(\frac{\partial f}{\partial x_j}
+
\frac{1}{\varepsilon}
\frac{\partial f}{\partial \xi_j}
\right)
\right|_{\tensor{\xi}=\frac{\mathbf{x}}{\varepsilon}}
=
\left.
\left(\frac{\partial f}{\partial x_j}
+
\frac{1}{\varepsilon}
f_{,j}
\right)
\right|_{\tensor{\xi}=\frac{\mathbf{x}}{\varepsilon}}
\end{equation}
and substituting the asymptotic expansion \eqref{eq:asymp_expansion} into the governing field equation \eqref{eq:field_equation}, one obtains
\begin{eqnarray}
\label{eq:RecursiveDifferentialProblemMicroDisplacement}
& &
\left.
\left\{
\varepsilon^{-2}
\left(
C_{ijkl}^m\, u_{k,l}^{(0)}
\right)_{,j}
+\varepsilon^{-1}
\left\{
\left[
C_{ijkl}^{m}
\left(
\frac{\partial u_{k}^{(0)}}{\partial x_l}+
u_{k,l}^{(1)}
\right)
\right]_{,j}+
\frac{\partial}{\partial x_j}
\left(
C_{ijkl}^{m}\,
u_{k,l}^{(0)}
\right)
\right\}
+
\right.
\right.
\nonumber\\
&&
\left.
\left.
+
\left[
C_{ijkl}^{m}
\left(
\frac{\partial u_{k}^{(1)}}{\partial x_l}+
u_{k,l}^{(2)}
\right)
\right]_{,j}
+
\frac{\partial}{\partial x_j}
\left[
C_{ijkl}^{m}
\left(
\frac{\partial u_{k}^{(0)}}{\partial x_l } + u_{k,l}^{(1)}
\right)
\right]
+
\varepsilon
\left\{
\left[
C_{ijkl}^m
\left(
\frac{\partial u_k^{(2)}}{\partial x_l}+
u_{k,l}^{(3)}
\right)
\right]_{,j} +
\right.
\right.
\right.
\nonumber\\
&&
\left.
\left.
\left.
 +
\frac{\partial}{\partial x_j}
\left[
C_{ijkl}^{m}
\left(
\frac{\partial u_{k}^{(1)}}{\partial x_l}
+
u_{k,l}^{(2)}
\right)
\right]
\right\}
+
O(\varepsilon^2)\right\}\right|_{\tensor{\xi}=\frac{\mathbf{x}}{\varepsilon}}
+b_i(\mathbf{x})=0.
\end{eqnarray}
A set of recursive partial differential problems can now be obtained from Eq.\eqref{eq:RecursiveDifferentialProblemMicroDisplacement} by equating terms at the different orders of $\varepsilon$.
Details about the solutions of such recursive differential problems can be found in \citep{Smyshlyaev2000,Bacigalupo2014,fantoni2017multi}. In particular, the solvability condition of the partial differential problem at the order $\varepsilon^{-2}$ in the class of $\mathcal{Q}$-periodic functions \citep{Bakhvalov1984}, implies that the solution at the order $\varepsilon^{-2}$ does not depend on the fast variable $\tensor{\xi}$, taking the form
\begin{equation}
\label{eq:micro_displ_zero_order}
u_k^{(0)}
\left(
\mathbf{x},
\tensor{\xi}
\right)
=
U_k
\left(
\mathbf{x}
\right)
\end{equation}
In force of solution \eqref{eq:micro_displ_zero_order} of the elastic differential problem at the order $\varepsilon^{-2}$, the solvability condition of resulting differential problem at the order $\varepsilon^{-1}$ and the $\mathcal{Q}$-periodicity of $\mathbb{C}^m$, lead to a solution at the first order of the form
\begin{equation}
\label{eq:micro_displ_order-1}
u_k^{(1)}
\left(
\mathbf{x},
\tensor{\xi}
\right)
=
N_{kpq_1}^{(1)}(\tensor{\xi})
\frac{\partial U_p\left(\mathbf{x}\right)}{\partial x_{q_1}}
\end{equation}
where $N_{kpq_1}^{(1)}$ is the first order perturbation function, which accounts for the influence of microstructural inhomogeneities.

Accounting for the solutions \eqref{eq:micro_displ_zero_order} and \eqref{eq:micro_displ_order-1} of differential problems at the order $\varepsilon^{-2}$ and $\varepsilon^{-1}$, respectively, and of the $\mathcal{Q}$-periodicity of the micro constitutive tensor components $C_{ijkl}^m$ and of the perturbation function $N_{kpq_1}^{(1)}$, at the order $\varepsilon^{2}$ the micro displacement solution takes the form
\begin{equation}
\label{eq:micro_displ_order0}
u_k^{(2)}
\left(
\mathbf{x},
\tensor{\xi}
\right)
=
N_{kpq_1q_2}^{(2)}(\tensor{\xi})
\frac{\partial^2 U_p\left(\mathbf{x}\right)}{\partial x_{q_1}\partial x_{q_2}}
\end{equation}
where $N_{kpq_1q_2}^{(2)}$ is the second order perturbation function which, again, depends only on the geometrical and physico-mechanical properties of the considered material.

Determination of perturbation functions derives from the solution of non-homogeneous recursive differential problems at the different order $\varepsilon$, being those problems referred to as {\itshape cell problems}.
Such differential problems are characterized by $\mathcal{Q}$-periodic source terms having a vanishing mean value over unit cell $\mathcal{Q}$.
Consequently, their solutions result to be $\mathcal{Q}$-periodic and they are enforced to have a zero mean value over $\mathcal{Q}$ in order to guarantee their uniqueness.

The following normalization condition is therefore imposed for all perturbation functions:
\begin{equation}
\label{eq:normalization_condition_pert_func}
\left\langle
\left(
\cdot
\right)
\right\rangle=
\frac{1}{|\mathcal{Q}|}
\int_{\mathcal{Q}}\left(\cdot\right) d\tensor{\xi}=0
\end{equation}
where $|\mathcal{Q}|$ represents the unit cell area, namely $|\mathcal{Q}|=\delta$.
The cell problem at the order $\varepsilon^{-1}$ in terms of the first order perturbation function $N_{kpq_1,l}^{(1)}$ reads
\begin{equation}
\label{eq:cell_probl_order-1}
\left(
C_{ijkl}^m N_{kpq_1,l}^{(1)}
\right)_{,j}
+
C_{ijpq_1,j}^m=0
\end{equation}
Once $N_{kpq_1}^{(1)}$ is determined, from Eq. (\ref{eq:RecursiveDifferentialProblemMicroDisplacement}) one obtains the cell problem at the order $\varepsilon^0$ which, symmetrized with respect to indices $q_1$ and $q_2$, takes the form
\begin{eqnarray}
\label{eq:cell_probl_order0}
&&
\left(C_{ijkl}^{m}\,N_{kpq_1q_2,l}^{(2)}
\right)_{,j}
+
\frac{1}{2}
\left[
\left(
C_{ijkq_2}^{m}\,
N_{kpq_1}^{(1)}
\right)_{,j}
+
C_{iq_1pq_2}^{m}
\right.
\nonumber\\
&&
\left.
+
C_{iq_2kl}^{m}\,
N_{kpq_1,l}^{(1)}
+
\left(
C_{ijkq_1}^{m}\,
N_{kpq_2}^{(1)}
\right)_{,j}
+
C_{iq_2pq_1}^{m}
+
C_{iq_1kl}^{m}\,
N_{kpq_2,l}^{(1)}
\right]=
\nonumber\\
&&
+
\frac{1}{2}
\left\langle
C_{iq_1pq_2}^m
+
C_{iq_2kl}^m\,
N_{kpq_1,l}^{(1)}
+
C_{iq_2pq_1}^{m}
+
C_{iq_1kl}^m\,
N_{kpq_2,l}^{(1)}
\right\rangle
\end{eqnarray}
Determination of perturbation functions as solutions of cell problems at the different orders of $\varepsilon$ allows obtaining the down-scaling relation, expressing the micro displacement field $\mathbf{u}(\mathbf{x},\tensor{\xi})$ as an asymptotic expansion of powers of $\varepsilon$ in terms of the $\mathcal{L}$-periodic macro field $\mathbf{U}(\mathbf{x})$ and its gradients. The down-scaling relation is expressed as
\begin{equation}
\label{eq:downscaling}
u_k(\mathbf{x},\tensor{\xi})
=
U_k(\mathbf{x})
+
\varepsilon
\,
N_{kpq_1}^{(1)}
(\tensor{\xi})
\frac{\partial U_p(\mathbf{x})}{\partial x_{q_1}}
+
\varepsilon^2
\,N_{kpq_1q_2}^{(2)}
(\tensor{\xi})
\frac{\partial^2 U_p(\mathbf{x})}{\partial x_{q_1}\partial x_{q_2}} + \mathcal{O}(\varepsilon^3)
\end{equation}
Denoting with $\tensor{\zeta}\in\mathcal{Q}$ a variable such that the vector $\varepsilon\tensor{\zeta}\in\mathcal{A}$ defines the translation of the medium with respect to $\mathcal{L}$-periodic body forces $\mathbf{b}(\mathbf{x})$ \citep{Smyshlyaev2000,Bacigalupo2014}, the
macro displacement field can in turn be expressed in terms of the micro field trough the up-scaling relation, as the mean value of $\mathbf{u}$ over unit cell $\mathcal{Q}$
\begin{equation}
\label{eq:upscaling}
U_k(\mathbf{x})
\doteq
\left\langle
u_k
\left(
\mathbf{x},
\frac{\mathbf{x}}{\varepsilon}
+
\tensor{\zeta}
\right)
\right\rangle_{\tensor{\zeta}}
\end{equation}
In particular, the variable $\tensor{\zeta}$ removes rapid fluctuations of coefficients and is such that
invariance property
\begin{equation}
\label{eq:invarance_property}
\left\langle
g(\tensor{\xi}+\tensor{\zeta})\right\rangle_{\tensor{\zeta}}
=
\frac{1}{|\mathcal{Q}|}
\int_{\mathcal{Q}} g(\tensor{\xi}+\tensor{\zeta})\,d\tensor{\zeta}
=
\frac{1}{|\mathcal{Q}|}
\int_{\mathcal{Q}} g(\tensor{\xi}+\tensor{\zeta})\,d\tensor{\xi}=
\left\langle g(\tensor{\xi})\right\rangle
\end{equation}
is proved to hold for $\mathcal{Q}$-periodic functions. Eq. \eqref{eq:invarance_property}, together with normalization condition \eqref{eq:normalization_condition_pert_func} enforced for all perturbation functions, leads to the up-scaling relation \eqref{eq:upscaling}.

\subsection{Closed form of the homogenized elastic tensor}
\label{SubSec:hom_constitutive_tensor}
%
The overall elastic  tensor is determined by means of a generalized macro-homogeneity condition, which establishes an energetic equivalence between the macro and the micro scales \citep{Smyshlyaev2000, Bacigalupo2014}.

The microscopic mean strain energy $\bar{E}_m$ is defined in terms of the micro strain energy density $\varphi_m$ as
\begin{equation}
\label{eq:micro_energy}
\bar{E}_m=\frac{1}{|\mathcal{Q}|}\int_{\mathcal{L}}
\int_{\mathcal{Q}}
\varphi_m
\left(
\mathbf{x},
\tensor{\xi}
\right)\,
d\tensor{\xi}
\,d\mathbf{x}=
\int_{\mathcal{L}}
\left\langle
\varphi_m\right\rangle\,d\mathbf{x}
\end{equation}
where $\varphi_m$ reads
\begin{equation}
\label{eq:micro_elastic_potential}
\varphi_m=
\frac{1}{2}C_{ijhk}^m
\frac{D u_i}{D x_j}
\frac{D u_h}{D x_k}=\frac{1}{2}C_{ijhk}^m\left(\delta_{ip}\delta_{jq_1}+N_{ipq_1,j}^{(1)}\right)
\left(\delta_{hs}\delta_{kr_1}+N_{hsr_1,k}^{(1)}\right)\frac{\partial U_p}{\partial x_{q_1}}
\frac{\partial U_s}{\partial x_{r_1}}+\mathcal{O}(\varepsilon)
\end{equation}
in view of the down-scaling relation \eqref{eq:downscaling}.
The strain energy at the macroscale is defined as
\begin{equation}
\label{eq:macro_energy}
E_M=\int_\mathcal{L}\varphi_M\,d\mathbf{x}
\end{equation}
in terms of the macro strain energy density $\varphi_M$ expressed in the form
\begin{equation}
\label{eq:macro_elastic_potential}
\varphi_M=
\frac{1}{2}C_{pq_1iq_2}
\frac{\partial U_p}{\partial x_{q_1}}
\frac{\partial U_i}{\partial x_{q_2}}.
\end{equation}
In Eq. \eqref{eq:macro_elastic_potential}, $C_{pq_1iq_2}$ denotes a generic component of the overall elastic tensor.
Denoting with $\bar{E}_{m}^0$ the mean micro strain energy truncated at the zeroth order, the generalized macro-homogeneity condition implies that
\begin{equation}
\label{eq:Hill_Mandell_gen}
\bar{E}_{m}^0\doteq E_M
\end{equation}
from which the components of the first order macroscopic tensor $C_{pq_1iq_2}$  can be defined in terms of the  components of the microscopic constitutive elastic tensor and in terms of the first order perturbation function. The closed form of the overall elastic tensor components is
\begin{equation}
\label{eq:Hom_elastic_tensor}
C_{pq_1iq_2}=
\left\langle
C_{rjkl}^{m}
\left(
N_{riq_2,j}^{(1)}+
\delta_{ir}\delta_{jq_2}
\right)
\left(
N_{kpq_1,l}^{(1)}+
\delta_{pk}\delta_{lq_1}
\right)
\right\rangle.
\end{equation}
Better approximation of the solution of the heterogeneous problem could be obtained by means of higher order homogenization techniques, modifyng the generalized macro-homogeneity condition \eqref{eq:Hill_Mandell_gen}  in a proper way.
Higher order approaches consider both higher order terms in the asymptotic expansion performed at the microscale and the non local character of the macroscopic constitutive tensor.
Alternatively, the solution of the average field equations of infinite order by means of perturbation methods allows to obtain a higher order approximation  by suitably truncating the asymptotic expansion of the macro field $\mathbf{U}$ in powers of $\varepsilon$ \citep{Bacigalupo2014}.
Higher order homogenization schemes could be explored in follow-up studies, motivating the development of higher order phase field formulations at the macroscale, for a consistent modelling of higher order effects across the material length scales.

\section{Damage evolution in a microstructured material trough a phase field model approach}
\label{Sec::phase_field}
%
Let  consider an arbitrary, homogeneous body $\Omega\in\mathbb{R}^2$. Note that for the sake of simplicity, we restrict our analysis to bidimensional geometries. As described in Sec. \ref{Sec::AsymHom}, vector $\mathbf{x}=x_1\mathbf{e}_1+x_2\mathbf{e}_2$ identifies the position of each material point inside the bulk at the macroscale. The boundary of the body $\Gamma\in\mathbb{R}$ is the union of Dirichlet $\Gamma_U$ and a Neumann $\Gamma_t$ parts, with $\Gamma_U\cup\Gamma_t=\Gamma$ and $\Gamma_U\cap\Gamma_t=\emptyset$.

Under the validity of the small strains assumption, the material response to the following quasi-static external actions is sought: body forces $\mathbf{b}(\mathbf{x})$ in $\Omega$, tractions $t$ on $\Gamma_t$, and displacements $\bar{U}$ on $\Gamma_U$.
Accordingly to the variational approach of fracture described in \citep{BourdinFrancfortMarigo,Miehe2010CMAME},
the total potential energy of a system $\Pi$ can be written as the balance between the internal contribution, which is the sum of a global energy storage functional and a dissipation functional due to damage, and the contribution due to external loading:
\begin{equation}
\label{eq:potential_energy}
\Pi \left(U_{i}, \mathfrak{d}\right) = \int_{\Omega} \varphi\left(H_{ij},\mathfrak{d}\right)\,\text{d}\Omega+\int_{\Omega}
\mathcal{G}_{c} \gamma \left(\mathfrak{d}, \frac{ \partial \mathfrak{d}}{\partial x_j}\right)\,\text{d}\Omega
-
\int_{\Gamma_t}t_i\,U_i\,\text{d}\Gamma
-
\int_{\Omega}b_i\,U_i\text{d}\Omega.
\end{equation}
Note that the first term on the right-hand side of Eq. \eqref{eq:potential_energy} is the global energy storage functional in which the strain energy density per unit volume $\varphi(H_{ij},\mathfrak{d})$ depends on both the macro displacement field trough $H_{ij}$ and the phase field $\mathfrak{d}$.
In the small-strain context, $H_{ij}$ is defined as the displacement gradient:
\begin{equation}
\label{eq:Macro_strain}
H_{ij}=
\frac{\partial U_i}{\partial x_j}.
\end{equation}

According to \citep{Miehe2010CMAME,MIEHE2010IJNME}, the phase field $\mathfrak{d}\in[0,1]$ is an internal auxiliary variable characterizing, for $\mathfrak{d}=0$, the undamaged condition and, for $\mathfrak{d}=1$, the fully damaged condition of the material.

The second term in Eq. (\ref{eq:potential_energy}) defines the work needed to create a diffusive damage and its rate provides the dissipation power due to damage in the bulk. In this context, $\mathcal{G}_C$ would correspond to a Griffith-type fracture energy, and $\gamma(\mathfrak{d},\partial \mathfrak{d}/\partial x_j)$ represents the crack surface density function depending on the phase field and on its gradient, which leads to the nonlocality of the formulation and therefore its mesh independency.
It can be expressed as
\begin{equation}
\label{eq:crack surface density}
\gamma
\left(
\mathfrak{d},
\frac{\partial\mathfrak{d}}{\partial x_j}
\right)=
\frac{1}{2\ell}
\mathfrak{d}^2
+
\frac{\ell}{2}
\left(
\frac{\partial\mathfrak{d}}{\partial x_j}
\right)^2
\end{equation}
where $\ell$ is an internal length scale parameter  governing the regularization/diffusion of damage \citep{kuhn2014,nguyen2016}. Here, we request to $\ell$ to be big enough to avoid crack localization and model situations where diffuse damage takes place.
In order to introduce damage only in tension, the energy density in the bulk is left undamaged in compression.
On the other hand, it is degraded in tension by assuming that only the Young's modulus of the matrix is affected by damage, i.e., by setting $E_m(\mathfrak{d})[(1-\mathfrak{d})^2+\mathcal{K}]E_{m,0}$, where $E_{m,0}$ is the undamaged value of the Young's modulus of the matrix and $\mathcal{K}$ is a residual material stiffness introduced to avoid numerical instabilities for $\mathfrak{d}=1$. The overall effective constitutive tensor components $C_{ijhk}(\mathfrak{d})$ is finally computed at the microscale by asymptotic homogenization and it is provided to the finite element computation scheme at the macroscale. In this regard, it is important to highlight that the elastic strain energy density is not degraded as in the standard phase field approach used in homogenous materials. Therefore, since the damaged constitutive tensor is the outcome of asymptotic homogenization, it will also depend upon the volumetric content and the shape of the inclusion/fiber.

The first order variation of functional \eqref{eq:potential_energy} reads:
\begin{eqnarray}
\label{eq:weak_form}
\delta\Pi(U_i,\delta U_i,\mathfrak{d},\delta\mathfrak{d}) &=&\int_{\Omega}H_{ij}\,C_{ijhk}\left(\mathfrak{d}\right)\delta H_{hk}\,\text{d}\Omega
+
\int_{\Omega}\frac{1}{2}H_{ij}\frac{\partial C_{ijhk}(\mathfrak{d})}{\partial \mathfrak{d}}H_{hk}\,\delta\mathfrak{d}
\,\text{d}\Omega
+\nonumber\\
&+&\int_{\Omega}
\frac{G_C}{\ell}
\mathfrak{d}\delta\mathfrak{d} \,\text{d}\Omega +
\int_{\Omega}
G_C\ell \frac{\partial \mathfrak{d}}{\partial x_j}\frac{\partial \delta\mathfrak{d}}{\partial x_j}
\,d\Omega - \int_{\Gamma_t}t_i\delta U_i\,d\Gamma - \int_{\Omega}b_i\,\delta U_i\,\text{d}\Omega
\end{eqnarray}
Integration by parts of the first and fourth terms of the right-hand side of Eq. (\ref{eq:weak_form}) allows to rephrase the variation $\delta \Pi$ as
\begin{eqnarray}
\label{eq:weak_form2}
&&\delta\Pi(U_i,\delta U_i,\mathfrak{d},\delta\mathfrak{d})=
\int_{\Gamma_t}
\frac{\partial U_i}{\partial x_j}\, C_{ijhk}\, n_k\, \delta U_h\, \text{d}\Gamma
-
\int_{\Omega} \frac{\partial}{\partial x_k}
\left(\frac{\partial U_i}{\partial x_j}\, C_{ijhk}
\right)
\delta U_h\,\text{d}\Omega
+
\nonumber\\
&&
\int_{\Omega}\frac{1}{2}H_{ij}\frac{\partial C_{ijhk}(\mathfrak{d})}{\partial \mathfrak{d}}H_{hk}\,\delta\mathfrak{d}
\,\text{d}\Omega
+
\int_{\Omega}
\frac{G_C}{\ell}
\mathfrak{d}\delta\mathfrak{d} \,\text{d}\Omega
+\int_{\Gamma}
G_C \,\ell\,\frac{\partial\mathfrak{d}}{\partial x_j}
n_j\,\delta\mathfrak{d}\,\text{d}\Gamma
+
\nonumber\\
&&
-
\int_{\Omega}
G_C\,\ell\frac{\partial^2\mathfrak{d}}{\partial x_j^2}\delta\mathfrak{d}\,\text{d}\Omega
- \int_{\Gamma_t}t_i\delta U_i\,d\Gamma - \int_{\Omega}b_i\,\delta U_i\,\text{d}\Omega
\end{eqnarray}
where symmetry properties of the elastic tensor $\mathbb{C}$ and the condition $\delta U_i =0 $ on $\Gamma_U$ have been exploited.

The Euler-Lagrange equations associated with the displacement and phase field problems yield the coupled field equations
\begin{eqnarray}
\label{eq:Euler-Lagrange equations}
&&\frac{\partial}{\partial x_j}
\left(
C_{ijhk}(\mathfrak{d})\,
\frac{\partial U_h}{\partial x_k}
\right) + b_i = 0,
\nonumber\\
&&
-\ell^2\frac{\partial^2\mathfrak{d}}{\partial x_j^2}+ \mathfrak{d}+\frac{\ell}{2 G_C}H_{ij} \frac{\partial C_{ijhk}(\mathfrak{d})}{\partial \mathfrak{d}} H_{hk}=0,
\end{eqnarray}
in the domain $\Omega$, along with Neumann boundary conditions:
\begin{eqnarray}
\label{eq:Neumann_boundary_conditions}
&& C_{ijhk}\frac{\partial U_h}{\partial x_k}n_j=t_i\hspace{0.5cm}on\hspace{0.2 cm}\Gamma_t\nonumber\\
&&
\frac{\partial\mathfrak{d}}{\partial x_j} n_j =0\hspace{0.5cm}on\hspace{0.2 cm}\Gamma
\end{eqnarray}
The solution of the nonlinear system \eqref{eq:Euler-Lagrange equations}, together with boundary conditions \eqref{eq:Neumann_boundary_conditions}, allows to describe the evolution of macroscopic mechanical quantities and of the phase field $\mathfrak{d}(\mathbf{x})$ inside the homogenized material.

\section{Numerical examples: influence of  the phase field internal length scale and  microstructure upon damage evolution}
\label{Sec::Benchmark}
\begin{figure}[h!]
  \centering
  \begin{tabular}{c c }
 \hspace{-0.75cm}
  \includegraphics[width=4cm]{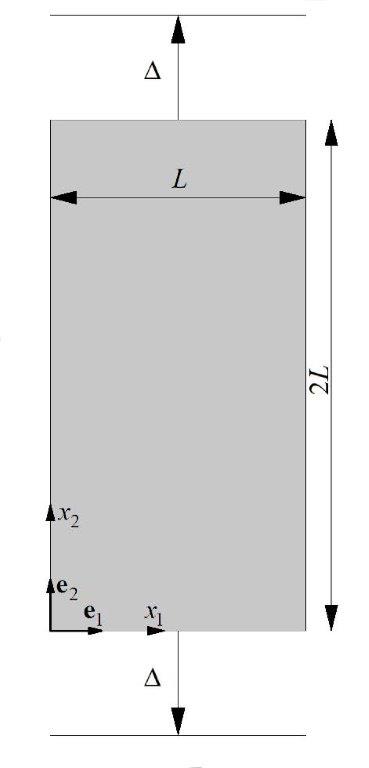}
  &
 \hspace{0.5cm}
  \includegraphics[width=4cm]{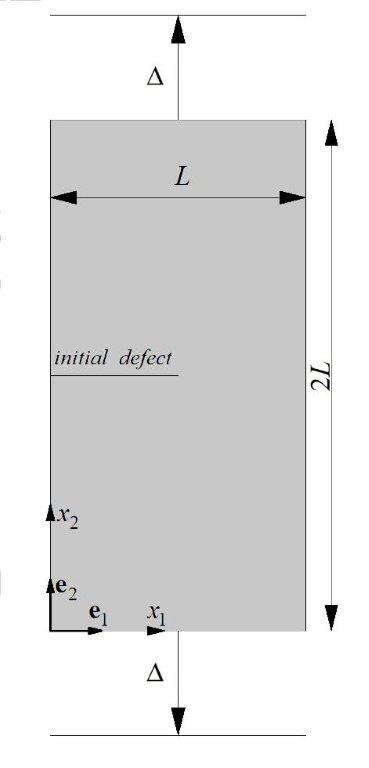}
  \\
   \hspace{-0.75cm} (a)  &  \hspace{0.5cm} (b)
  \end{tabular}
  \caption{\it  Specimens in plane strain conditions having width equal to $L=1$ mm and height equal to $2L$ subjected to imposed displacement $\Delta$ along the upper and lower boundaries. (a) Plane specimen, (b) specimen with an initial edge crack. }
  \label{Fig::Specimen}
\end{figure}

%
%
The main objective of the present section is to show the dependence of the evolution of damage inside the material to the shape and volumetric content of the inclusion, as well as to the phase field internal length $\ell$.
To this aim, let consider the specimens sketched in Fig. \ref{Fig::Specimen}, having width $L=1$ mm and height $2L$, both subjected to uniform tensile loading through the application of an imposed far-field displacement $\Delta$ to all the finite element nodes belonging to the lower and upper sides.
The case depicted in Fig. \ref{Fig::Specimen}-(b) distinguishes from the one of Fig. \ref{Fig::Specimen}-(a)  because of the presence of an initial defect, which is a horizontal straight edge crack with mouth at $x_2=L$ and crack tip at $x_1=L/2$.
In both cases, the specimen is supposed to be constituted by a microstructured material composed by a matrix embedding a circular or a square inclusion. Damage is assumed to develop only within the matrix.

In order to accelerate coupling between the macroscale and the microscale computations, asymptotic homogenization of the microstructured material as described in Section \ref{Sec::AsymHom} is performed off-line, for a set of admissible values of the phase field variable $\mathfrak{d}$ ranging from zero to unity.
The numerical results of the coupled problem \eqref{eq:Euler-Lagrange equations} for the test problem without initial defect and with an initial defect are discussed in Sections \ref{SubSec::Case1} and \ref{SubSec::Case2}, respectively.

\subsection{Homogenization of the material with microstructure}
\begin{figure}[h!]
  \centering
  \begin{tabular}{c c }
 \hspace{-0.75cm}
  \includegraphics[width=4cm]{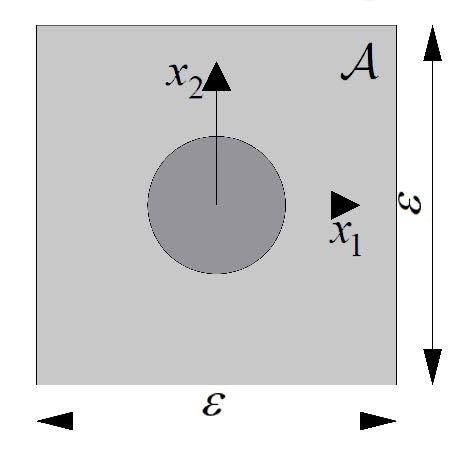}
  &
 \hspace{0.5cm}
  \includegraphics[width=4cm]{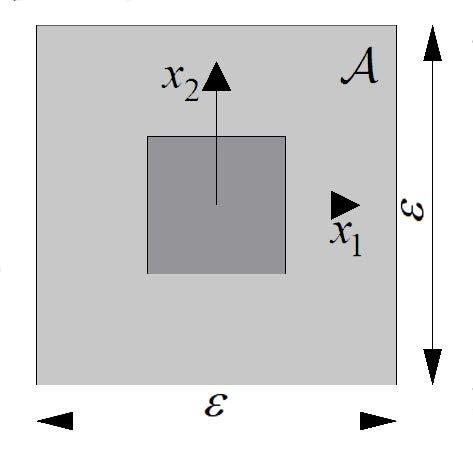}
  \\
   \hspace{-0.75cm} (a)  &  \hspace{0.5cm} (b)
  \end{tabular}
  \caption{\it Periodic cell $\mathcal{A}$ of characteristic size $\varepsilon$ made of an Aluminum matrix with a Silicon carbide inclusion. (a) Inclusion with circular shape, (b) Inclusion with a square shape.}
  \label{Fig::PeriodicCell}
\end{figure}

Specimens depicted in Fig. \ref{Fig::Specimen} are supposed to be made by a material with microstructure whose periodic cell is shown in Fig. \ref{Fig::PeriodicCell}.
The periodic cell is composed of an Aluminum matrix undergoing damage, with an undamaged linear elastic constitutive behavior characterized by a Young's modulus $E_{m,0}=E_{Al}=60$ GPa and Poisson ratio $\nu_{Al}=0.3$, and a linear elastic inclusion corresponding to Silicon carbide with Young's modulus $E_{SiC}=340$ GPa and $\nu_{SiC}=0.18$.

Under plane strain conditions, the elastic tensor of the undamaged matrix $\mathbb{C}_{Al}$ and of the inclusion $\mathbb{C}_{SiC}$ are:
\begin{equation}
\label{eq:elastic_tensors_materials}
\mathbb{C}_{Al}=
\left(
\begin{array}{c c c}
8.077  & 3.461 & 0 \\
3.461  & 8.077 & 0\\
0 & 0 & 2.308
\end{array}
\right) \times 10^4\,\text{MPa},
\hspace{0.5 cm}
\mathbb{C}_{SiC}=
\left(
\begin{array}{c c c}
36.917  & 8.104 & 0 \\
8.104  & 36.917 & 0\\
0 & 0 & 14.407
\end{array}
\right) \times 10^4\,\text{MPa}
\end{equation}

Two different inclusion shapes have been taken into account, namely circular (Fig. \ref{Fig::PeriodicCell}-(a)) with five different values of the volumetric content $f=1/4,1/8,1/16,1/32,1/100$, and square (Fig. \ref{Fig::PeriodicCell}-(b)) with a volumetric content $f=1/4$.
The volumetric content represents the ratio between the area of the inclusion and that of the periodic cell $\mathcal{A}$ in the cross-section of the specimen, as customary.

After computing the perturbation functions $N_{kpq_1}^{(1)}$ as the solutions of the cell problem \eqref{eq:cell_probl_order-1}, the overall elastic constitutive tensor has been computed by means of the closed form \eqref{eq:Hom_elastic_tensor}. Considering, for example, a volumetric content $f=1/4$, we get:
\begin{equation}
\label{eq:Hom_elastic_tensor_benchmark}
\mathbb{C}_{circ}=
\left(
\begin{array}{c c c}
10.896 & 4.104 & 0 \\
4.104 & 10.896 & 0 \\
0 & 0 & 3.148
\end{array}
\right)\times 10^4 \text{MPa},
\hspace{0.5cm}
\mathbb{C}_{sq}=
\left(
\begin{array}{c c c}
10.876 & 4.014 & 0 \\
4.0139 & 10.876 & 0\\
0 & 0 & 3.0464
\end{array}
\right)\times 10^4 \text{MPa}
\end{equation}
where subscripts $_{circ}$ and $_{sq}$ stand for circular and square topology, respectively.

\begin{figure}[h!]
  \centering
  \begin{tabular}{c c }
 \hspace{-0.75cm}
  \includegraphics[width=8cm]{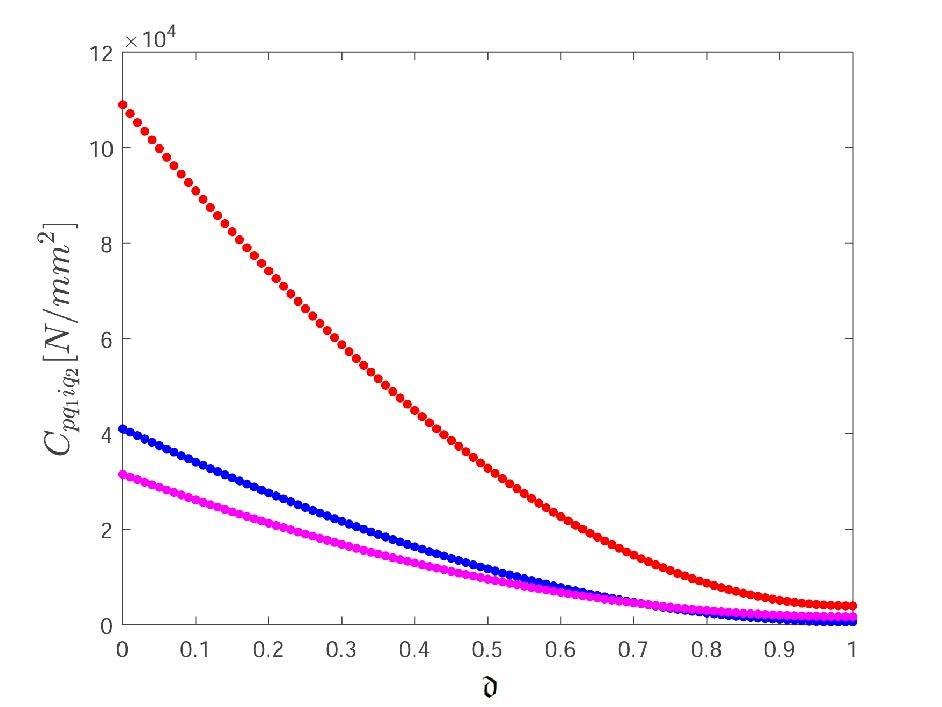}
  &
 \hspace{0.5cm}
  \includegraphics[width=8cm]{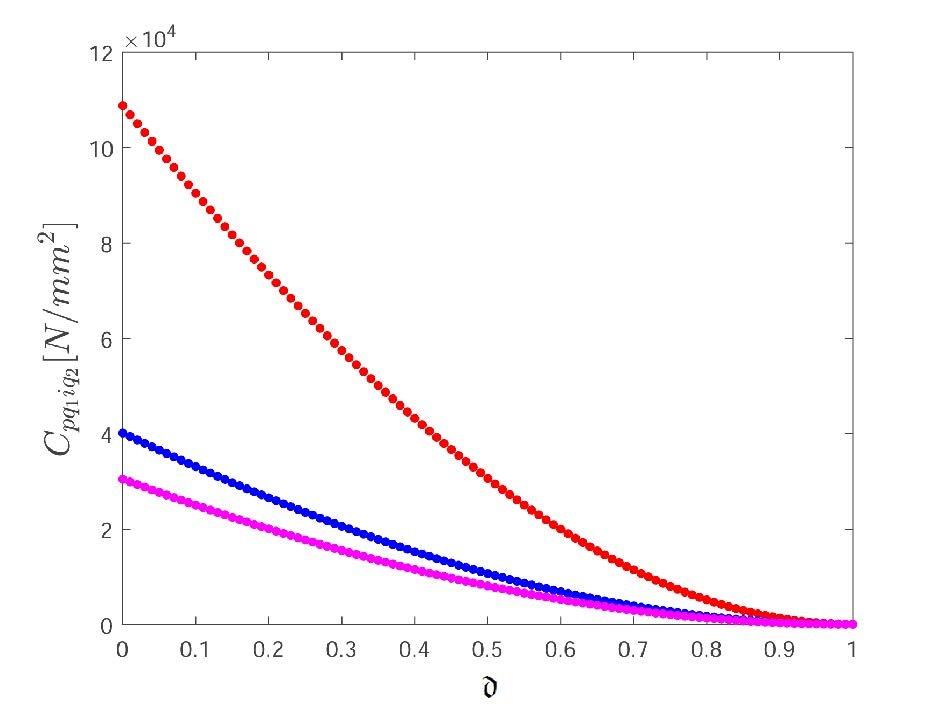}
  \\
   \hspace{-0.75cm} (a)  &  \hspace{0.5cm} (b)
  \end{tabular}
  \caption{\it Components of the overall elastic constitutive tensor $\mathbb{C}$: $C_{1111}$ (red), $C_{1122}$ (blue), and $C_{1212}$ (magenta). (a) Circular inclusion. (b) Square inclusion.}
  \label{Fig::GlobConst_vs_esse}
\end{figure}

Adopting a residual stiffness $\mathcal{K}=0.005$ in correspondence of $\mathfrak{d}=1$, the elastic constitutive tensor of the damaged matrix $\mathbb{C}_{Al}$ have been multiplied by degradation function $g(\mathfrak{d})=(1-\mathfrak{d})^2+\mathcal{K}$, with the phase field variable in the range $0\leq\mathfrak{d}\leq 1$.
The components of the overall elastic tensor $\mathbb{C}$ have then been computed for each value of $\mathfrak{d}$.
Fig. \ref{Fig::GlobConst_vs_esse} shows, for example, the components of $\mathbb{C}_{circ}$ (Fig. \ref{Fig::GlobConst_vs_esse}-(a)) and of $\mathbb{C}_{sq}$ (Fig. \ref{Fig::GlobConst_vs_esse}-(b)) as a function of the phase field $\mathfrak{d}$, for a microstructure configuration characterized by a volumetric content $f=1/4$.
In correspondence of $\mathfrak{d}=0$, the overall elastic constitutive tensors $\mathbb{C}_{circ}$ and $\mathbb{C}_{sq}$  correspond to those of the undamaged material, given by  Eq. \eqref{eq:Hom_elastic_tensor_benchmark}. Their values monotonically decrease as $\mathfrak{d}$ approaches unity, where the value $\mathcal{K}$ is retrieved.

Considering a constant value for $G_C=6$ N/mm for the homogenized material, the coupled system of equations \eqref{eq:Euler-Lagrange equations} has been numerically solved for each value of the imposed displacement $\Delta$, for the unnotched and the notched specimen configurations.
Details regarding the finite element framework adopted to obtain the numerical solution, and implemented in the finite element software FEAP \citep{Zienkiewicz1977}, are described in the Appendix \ref{Appendix:finite element}.
Five different values of the internal length scale $\ell$ have been considered, in order to investigate the role of this parameter in triggering damage inside the microstructured material. To avoid strain localization which would violate the applicability of homogenization, values of $\ell$ have been set intermediate between the value of the macroscopic length scale $L$ and the value of the microscopic one $\varepsilon$, namely $\ell=0.05,0.1,0.2,0.4$ mm.
\subsection{ Case 1: tensile test on a plain unnotched specimen}
\label{SubSec::Case1}
\begin{figure}[h!]
  \centering
  \begin{tabular}{c c}
 \hspace{-0.75cm}
  \includegraphics[width=8cm]{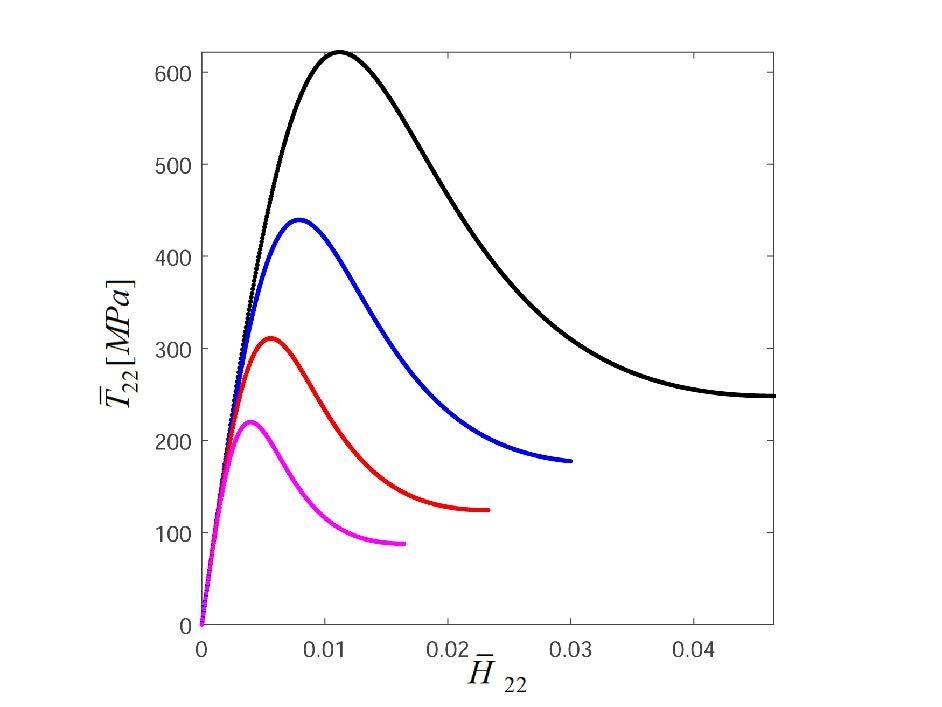}
  &
 \hspace{0.5cm}
  \includegraphics[width=8cm]{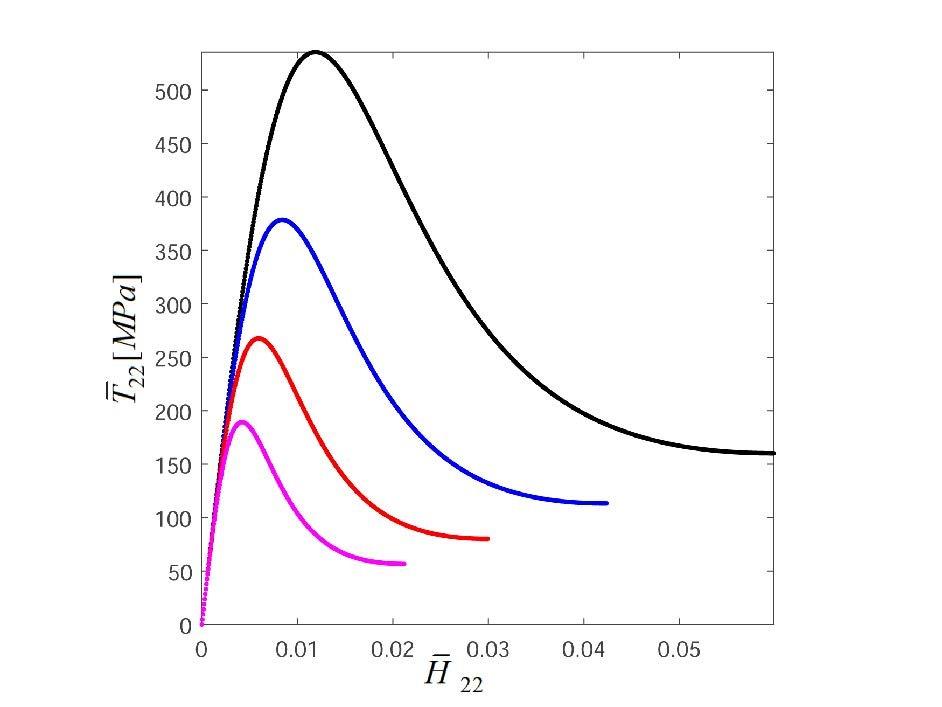}
  \\
   \hspace{-0.75cm} (a)  &  \hspace{0.5cm} (b)
   \\
 \hspace{-0.75cm}
  \includegraphics[width=8cm]{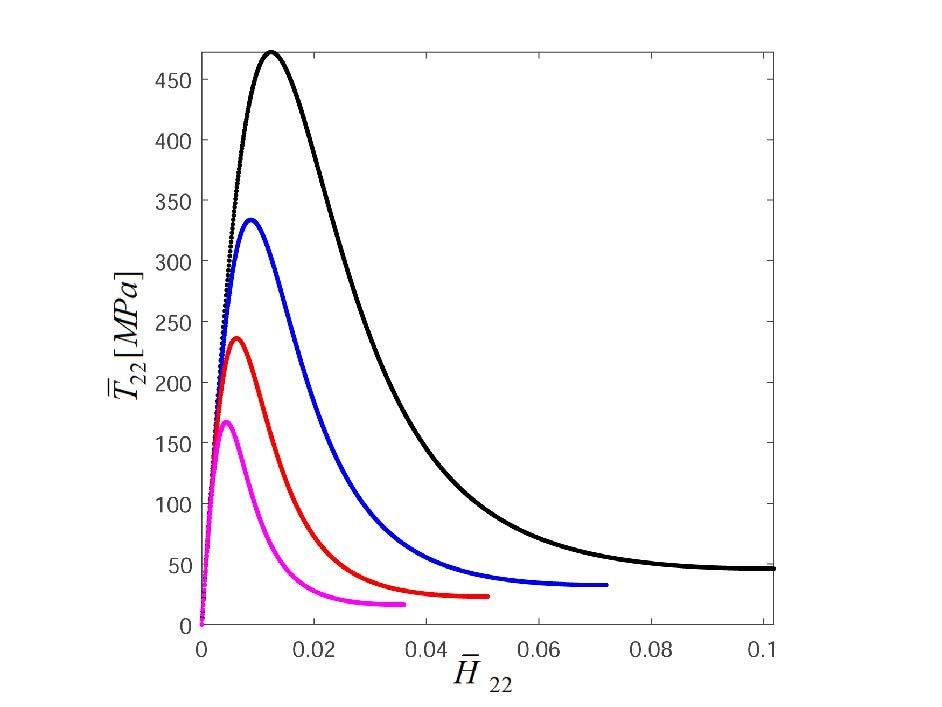}
  &
 \hspace{0.5cm}
  \includegraphics[width=8cm]{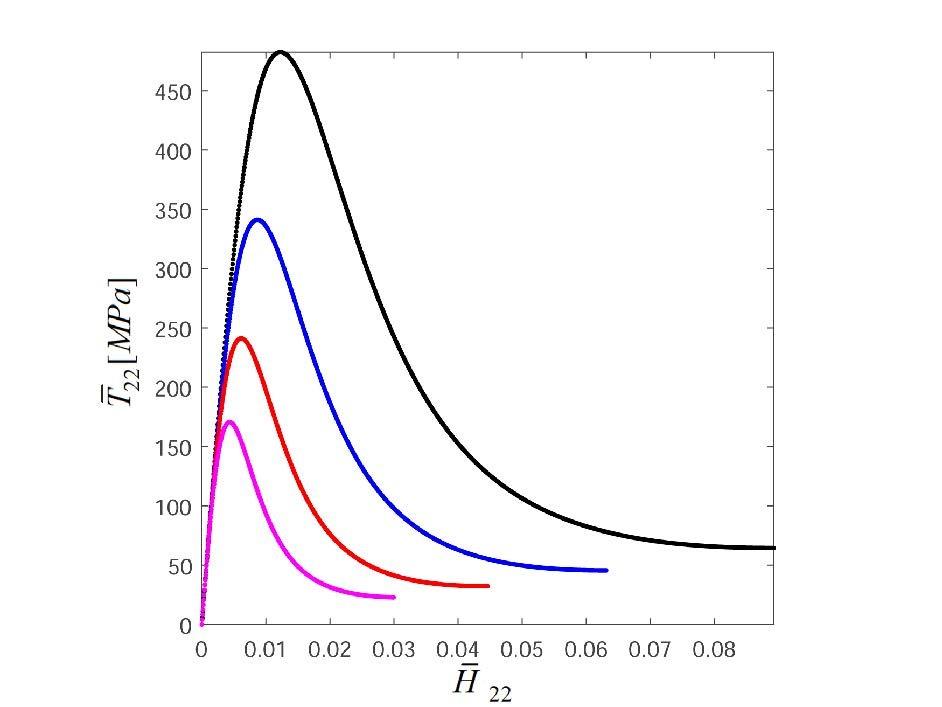}
  \\
   \hspace{-0.75cm} (c)  &  \hspace{0.5cm} (d)
   \\
 \hspace{-0.75cm}
  \includegraphics[width=8cm]{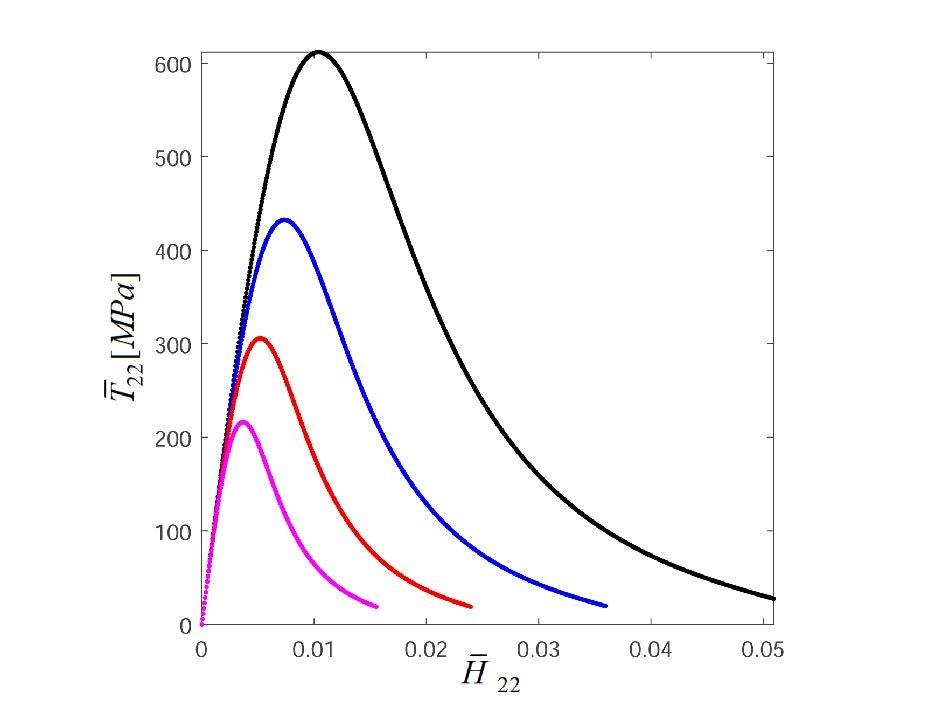}
  &
 \hspace{0.5cm}
  \includegraphics[width=8cm]{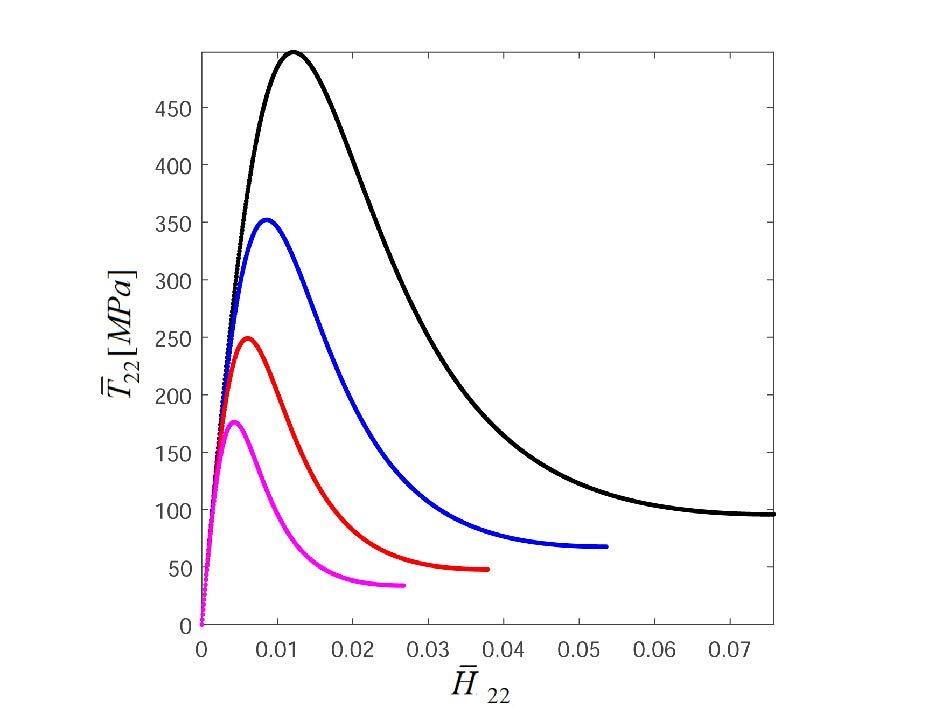}
  \\
   \hspace{-0.75cm} (e)  &  \hspace{0.5cm} (f)
   \\
  \end{tabular}
  \caption{\it Average stress $\bar{T}_{22}$ vs average deformation $\bar{H}_{22}$ for different values of the internal length scale $\ell$: $\ell=0.05$ (black), $\ell=0.1$ (blue), $\ell=0.2$ (red), and $\ell=0.4$ (magenta). (a) circular inclusion and $f=1/4$, (b) circular inclusion and $f=1/8$, (c) circular inclusion and $f=1/16$, (d) circular inclusion and $f=1/32$, (e) circular inclusion and $f=1/100$, (f) square inclusion and $f=1/4$. }
  \label{Fig::NO_int_diff_lo}
\end{figure}

One considers the plain unnotched specimen depicted in Fig. \ref{Fig::Specimen}-(a), where a uniaxial tensile test is performed.
Fig. \ref{Fig::NO_int_diff_lo} depicts the variation of the homogenized mean stress $\bar{T}_{22}$ with respect to the homogenized mean strain $\bar{H}_{22}$, by varying the phase field length scale $\ell$, for a circular inclusion with a different volumetric content $f$ in the subfigures from (a) to (e), and for a the square inclusion with $f=1/4$ in the subfigure (f).
In particular, the average stress $\bar{T}_{22}$ is computed as the sum of the reaction forces on the upper (or lower) boundary of the specimen, divided by the width $L$ of the specimen and its unit out-of-plane thickness. The average strain $\bar{H}_{22}$ is computed by dividing the imposed far-field displacement $\Delta$ by $L$.

According to previous results reported in the literature \citep{kuhn2014,nguyen2016,paggi2017}, the apparent strength $\bar{T}_{22_{max}}$, representing the maximum of the average stress-strain curves, increases as the length scale $\ell$ decreases. This trend can be intuitively explained by the fact that larger values of $\ell$ indicate a more widespread damage in the bulk. However, by comparing the different subfigures in Fig. \ref{Fig::NO_int_diff_lo}, one notices that the present model provides an apparent strength also dependent upon the volumetric content of the reinforcement and on the shape of the inclusion, which is a major novelty as compared to the standard phase field approach for homogeneous materials.
Examining each stress-strain curve more in detail, the system response is almost linear till $\bar{T}_{22_{max}}$ is reached. The stiffness of such an elastic part is the same by varying $\ell$, which is consistent with the fact that all the curves in a subfigure present the same microstructure with the same volumetric content $f$.
For an axial strain larger than the value of the strain corresponding to the apparent strength, the post-peak response presents softening. This is a major difference from the prediction of the phase field in the case of a homogenous material, and it is caused by the progressive stress transfer from the damaged matrix to the undamaged inclusion, by increasing damage. This redistribution of stresses cannot take place if the material is homogenous.
The value of the residual stress are increasing by reducing the length scale $\ell$, as an effect of the different evolution and spread of damage.
By comparing Figs. \ref{Fig::NO_int_diff_lo}-(a) and \ref{Fig::NO_int_diff_lo}-(f), which refer to the same value of $f$, but correspond to different inclusion shapes, one can notice how the apparent strength and the corresponding value of $\bar{H}_{22_{max}}$  are slightly influenced by the change of topology, while the post peak behavior is greatly affected by the inclusion shape, showing steeper branches and lower residual stresses in the case of a square inclusion.
Furthermore, Fig. \ref{Fig::NO_int_diff_lo}-(f) shows that the residual stresses are slightly influenced by $\ell$, for a square topology of the cell inclusion.

\begin{figure}[h!]
  \centering
  \begin{tabular}{c c}
 \hspace{-0.75cm}
  \includegraphics[width=8cm]{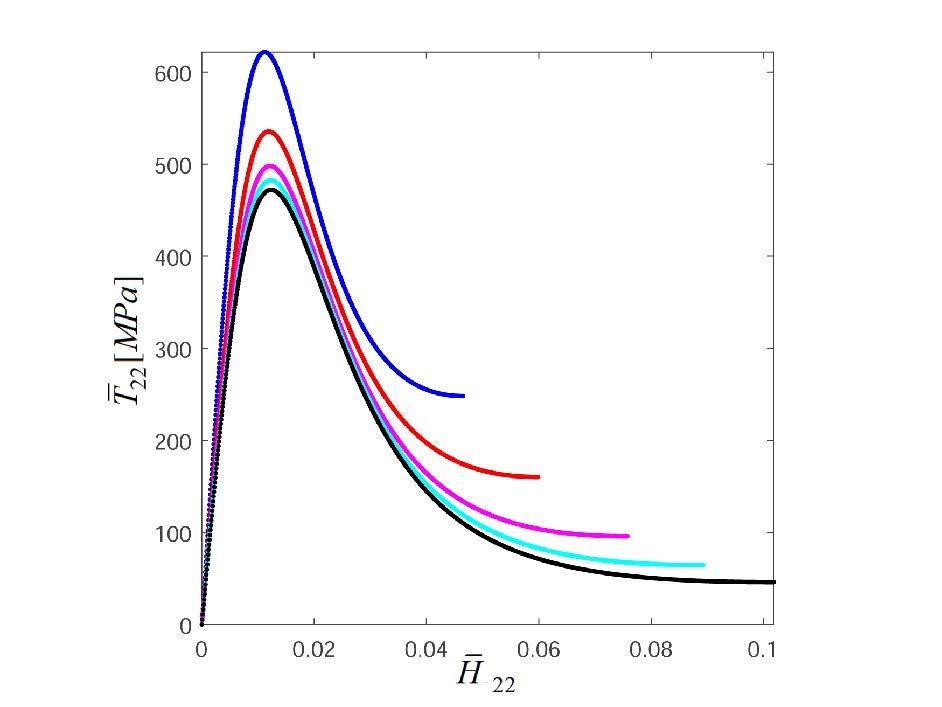}
  &
 \hspace{0.5cm}
  \includegraphics[width=8cm]{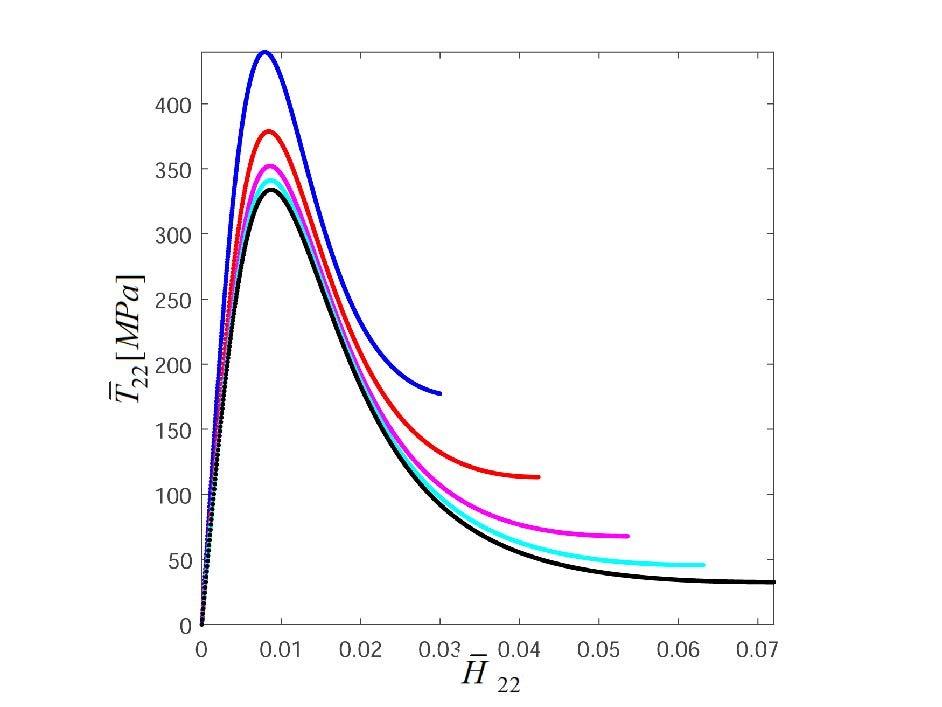}
  \\
   \hspace{-0.75cm} (a)  &  \hspace{0.5cm} (b)
   \\
 \hspace{-0.75cm}
  \includegraphics[width=8cm]{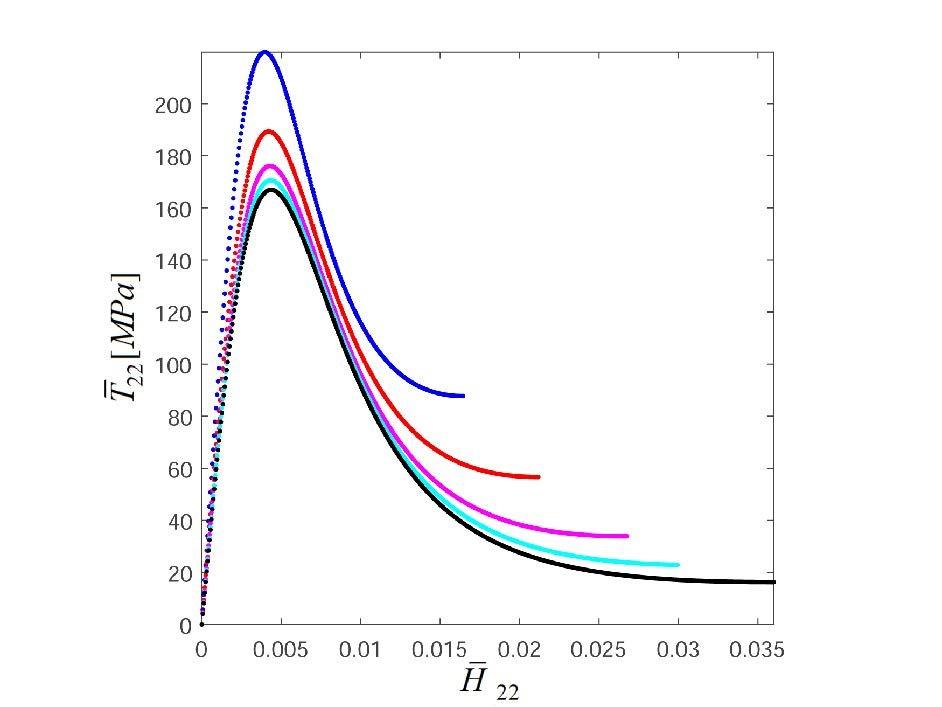}
  &
 \hspace{0.5cm}
  \includegraphics[width=8cm]{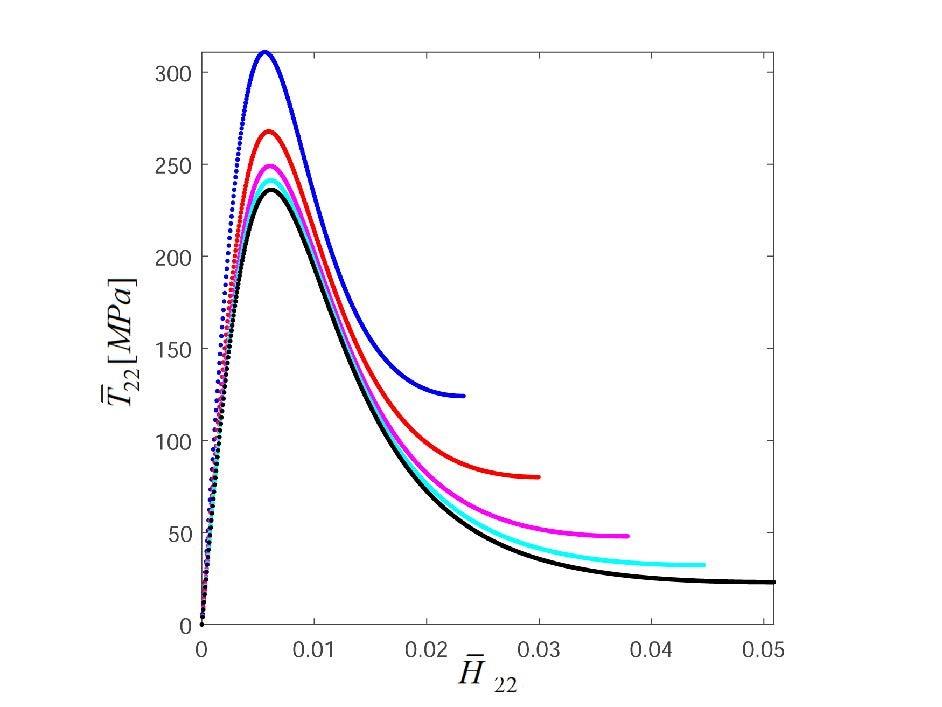}
  \\
   \hspace{-0.75cm} (c)  &  \hspace{0.5cm} (d)
  \end{tabular}
  \caption{\it Average stress $\bar{T}_{22}$ vs average deformation $\bar{H}_{22}$ for circular inclusion at different volume fractions:
  $f=1/4$ blue curve, $f=1/8$ red curve, $f=1/16$ magenta curve, $f=1/32$ cyan curve, $f=1/100$ black curve. (a) $\ell=0.05$ mm,(b) $\ell=0.1$ mm, (c) $\ell=0.2$ mm, (d) $\ell=0.4$ mm.
 }
  \label{Fig::NO_int_diff_f}
\end{figure}

\begin{figure}[h!]
  \centering
  \begin{tabular}{c c}
 \hspace{-0.75cm}
  \includegraphics[width=8cm]{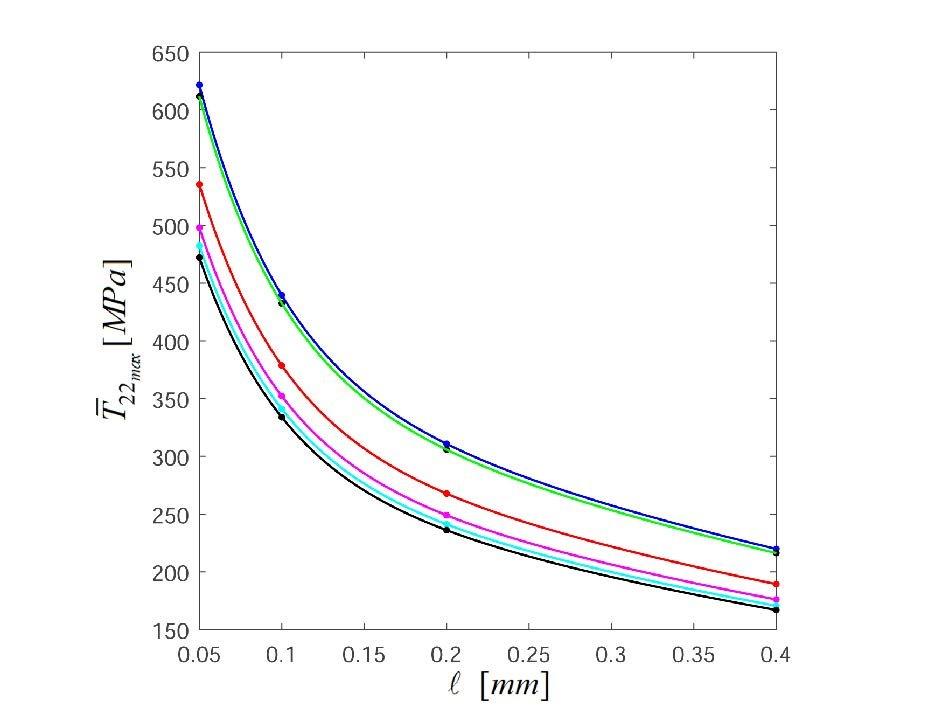}
  &
 \hspace{0.5cm}
  \includegraphics[width=8cm]{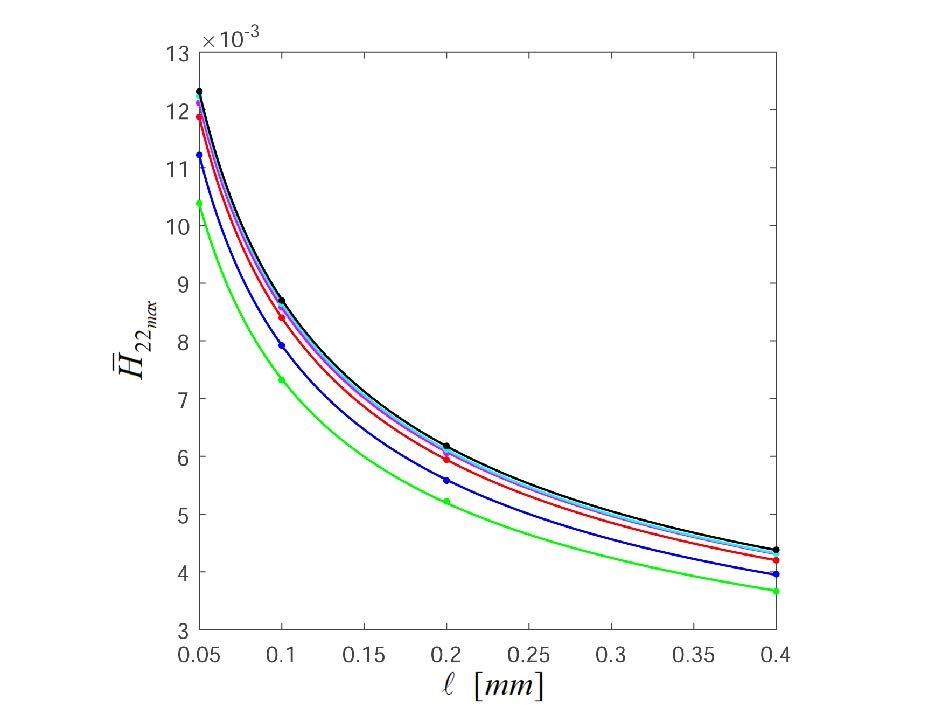}
  \\
   \hspace{-0.75cm} (a)  &  \hspace{0.5cm} (b)
  \end{tabular}
  \caption{\it Apparent strength $\bar{T}_{22_{max}}$ vs. $\ell$ (a), and the corresponding average deformation $\bar{H}_{22_{max}}$ vs. $\ell$ (b). Circular inclusion and $f=1/4$ (blue), $f=1/8$ (red), $f=1/16$ (magenta), $f=1/32$ (cyan), $f=1/100$ (black), square inclusion and $f=1/4$ (green).}
  \label{Fig::NO_int_Sigma_Eps_max}
\end{figure}

To highlight the effect of the volumetric content, Fig. \ref{Fig::NO_int_diff_f} depicts, for a circular inclusion, the average stress $\bar{T}_{22}$ as a function of $\bar{H}_{22}$ by keeping $\ell$ constant and varying $f$.
The apparent strength  $\bar{T}_{22_{max}}$  increases as volume fraction $f$ increases, and this is coherent with the choice to weaken the constitutive properties only for the matrix and not for the inclusion.
For the same reason, residual stresses are greater as volume fraction $f$ increases.
To summarize the above results, Fig. \ref{Fig::NO_int_Sigma_Eps_max} shows how the apparent strength of the system, $\bar{T}_{22_{max}}$, and the correspondent values of the strain, $\bar{H}_{22_{max}}$, depend on $\ell$, for the different inclusion volumetric contents and also for different shape of the inclusion.
Curves referring to the square topology and $f=1/4$ show   lower values than those for a circular inclusion and the same value of $f$, both for $\bar{T}_{22_{max}}$ and $\bar{H}_{22_{max}}$.

\subsection{ Case 2: tensile test for a specimen with an edge notch}
\label{SubSec::Case2}
\begin{figure}[h!]
  \centering
  \begin{tabular}{c c}
 \hspace{-0.75cm}
  \includegraphics[width=8cm]{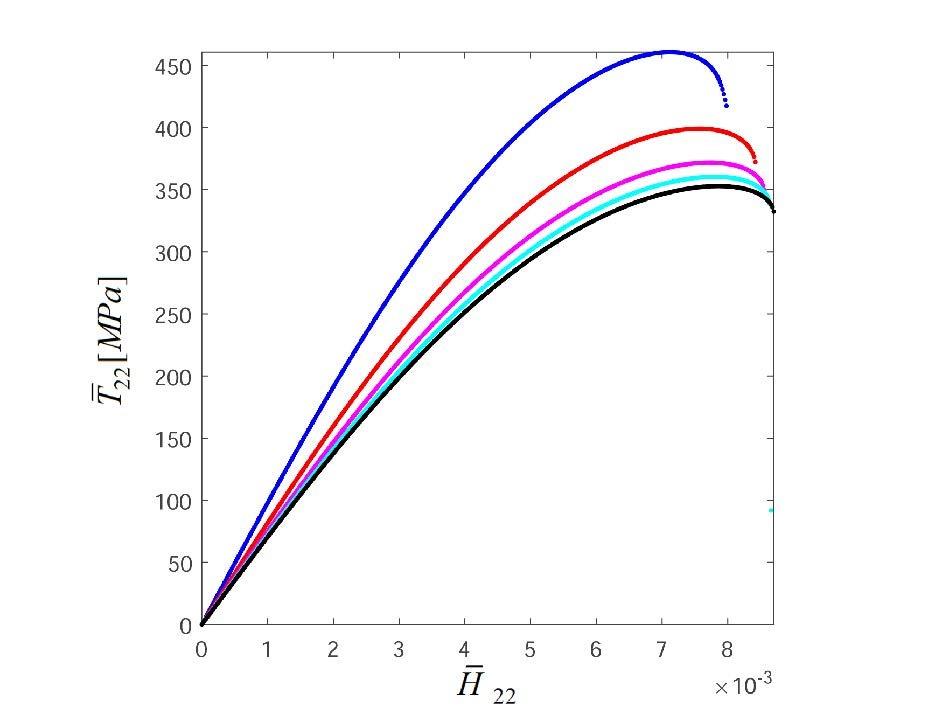}
  &
 \hspace{0.5cm}
  \includegraphics[width=8cm]{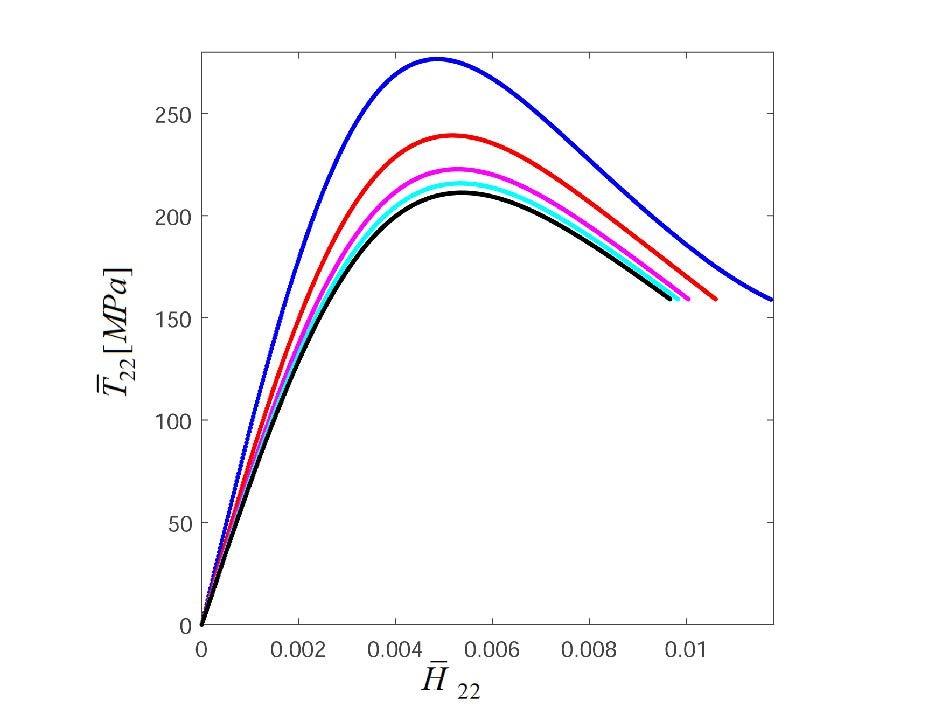}
  \\
   \hspace{-0.75cm} (a)  &  \hspace{0.5cm} (b)
  \end{tabular}
  \caption{\it Average stress $\bar{T}_{22}$ vs. average strain $\bar{H}_{22}$ for a material with circular inclusion and different volumetric content: $f=1/4$ (blue), $f=1/8$ (red), $f=1/16$ (magenta), $f=1/32$ (cyan), $f=1/100$ (black). (a) $\ell=0.1$ mm (b) $\ell=0.4$ mm. }
  \label{Fig::SI_int_diff_f}
\end{figure}

\begin{figure}[h!]
  \centering
  \begin{tabular}{c c}
 \hspace{-0.75cm}
  \includegraphics[width=8cm]{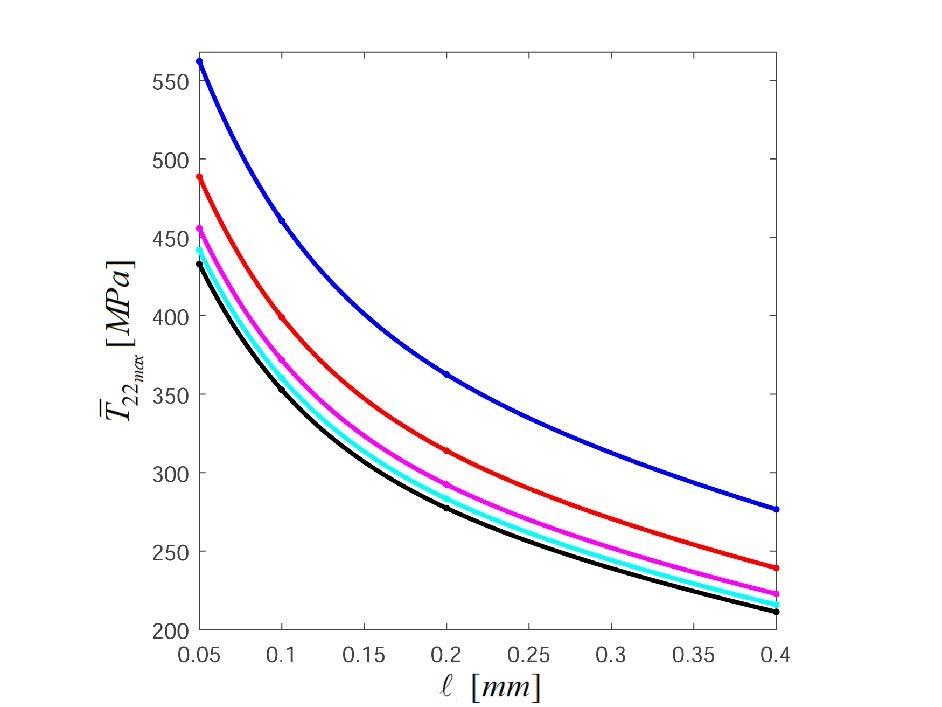}
  &
 \hspace{0.5cm}
  \includegraphics[width=8cm]{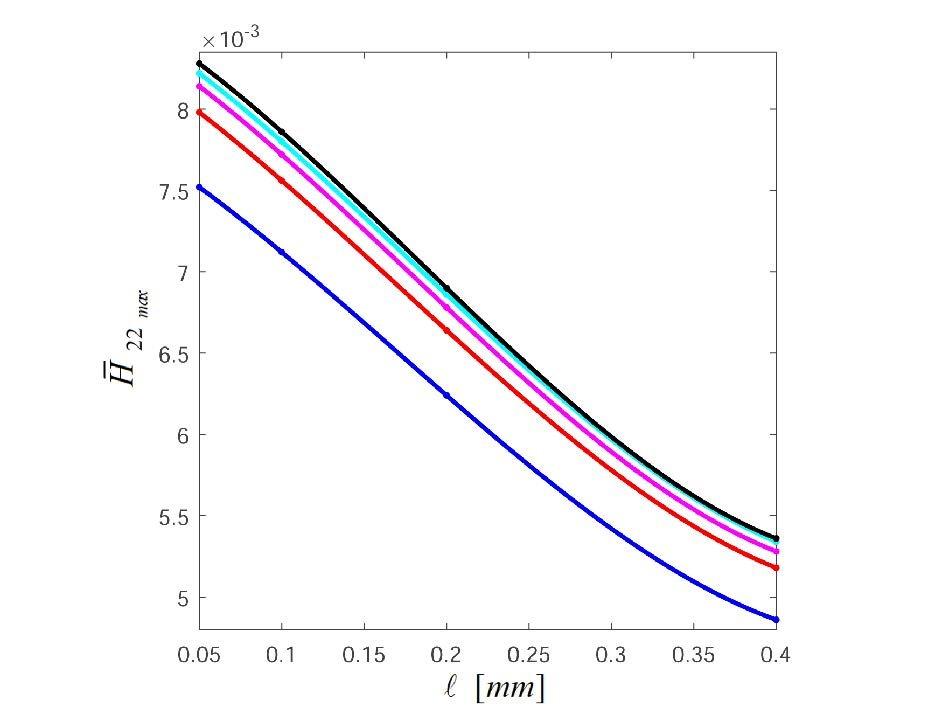}
  \\
   \hspace{-0.75cm} (a)  &  \hspace{0.5cm} (b)
  \end{tabular}
  \caption{\it Apparent strength $\bar{T}_{22_{max}}$ vs. $\ell$ (a), and corresponding values of $\bar{H}_{22_{max}}$ vs. $\ell$ (b) for a material microstructure with a circular inclusion and $f=1/4$ (blue), $f=1/8$ (red), $f=1/16$ (magenta), $f=1/32$ (cyan), $f=1/100$ (black).}
  \label{Fig::SI_int_Sigma_Eps_max}
\end{figure}

Let consider a uniaxial tensile test for the notched specimen depicted in Fig. \ref{Fig::Specimen}-(b).
Fig. \ref{Fig::SI_int_diff_f} shows the average stress $\bar{T}_{22}$ as a function of average strain $\bar{H}_{22}$ for this test, where $\bar{T}_{22}$ and $\bar{H}_{22}$ are computed as for the specimen without an initial notch.
In particular, Fig. \ref{Fig::SI_int_diff_f} shows the homogenized stress-strain curves for a material microstructure with a circular inclusion and different values of $f$, for two distinct values of the internal length scale $\ell$, namely $\ell=0.1$ mm (Fig. \ref{Fig::SI_int_diff_f}-(a)) and $\ell=0.4$ mm (Fig. \ref{Fig::SI_int_diff_f}-(b)).
Once again, the two-scale coupled formulation predicts an apparent strength $\bar{T}_{22_{max}}$ which is an increasing function of $f$, while the correspondent average deformation $\bar{H}_{22_{max}}$ is a decreasing function of $f$.
The mechanical response is quite affected by $\ell$ as can be assessed by comparing Figs. \ref{Fig::SI_int_diff_f}-(a) and \ref{Fig::SI_int_diff_f}-(b). Softening branches occur for $\ell=0.4$ mm, while for $\ell=0.1$ mm the post-peak response is quite brittle.
Fig. \ref{Fig::SI_int_Sigma_Eps_max} shows $\bar{T}_{22_{max}}$ and the corresponding strain $\bar{H}_{22_{max}}$ as functions of $\ell$, for all the different values of $f$ herein considered.
The apparent strength of the system and the corresponding  homogenized strain are monotonically decreasing functions of the internal length scale $\ell$. Moreover, the comparison between Figs. \ref{Fig::NO_int_Sigma_Eps_max} and \ref{Fig::SI_int_Sigma_Eps_max} shows that, as expected, $\bar{T}_{22_{max}}$ and $\bar{H}_{22_{max}}$ are lower than the corresponding values for the unnotched specimen.

\begin{figure}[h!]
  \centering
  \begin{tabular}{c c c }
 \hspace{-0.75cm}
  \includegraphics[width=3.1cm]{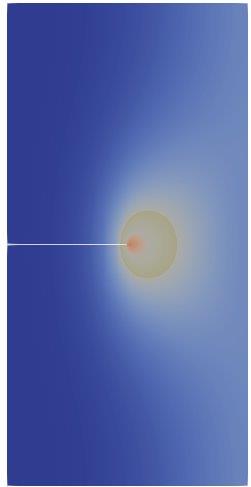}
  &
 \hspace{0.5cm}
  \includegraphics[width=3.15 cm]{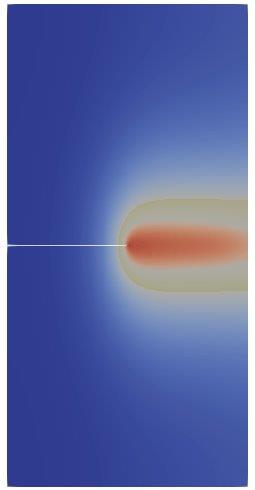}
  &
 \hspace{-0.4cm}
  \includegraphics[width=1.4cm]{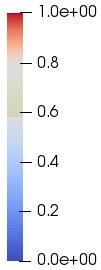}
  \\
   \hspace{-0.75cm} (a)  &  \hspace{0.5cm} (b) &
  \end{tabular}
  \caption{\it  Contour plot of $\mathfrak{d}$ for $\ell=0.1$ mm, circular inclusion and $f=1/4$. (a)  $\bar{H}_{22}=0.0062 $  (b) $\bar{H}_{22}=0.0072 $  }
  \label{Fig::Specimen_l0_01}
\end{figure}

\begin{figure}[h!]
  \centering
  \begin{tabular}{c c c }
 \hspace{-0.75cm}
  \includegraphics[width=3.04cm]{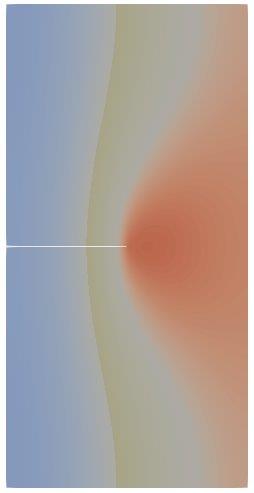}
  &
 \hspace{0.5cm}
  \includegraphics[width=3 cm]{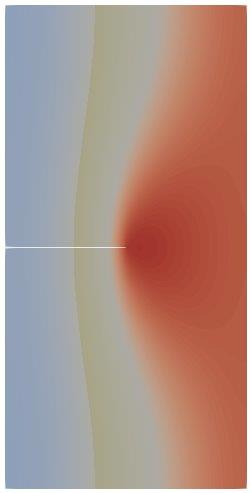}
  &
 \hspace{-0.4cm}
  \includegraphics[width=1.4cm]{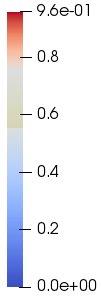}
  \\
   \hspace{-0.75cm} (a)  &  \hspace{0.5cm} (b) &
  \end{tabular}
  \caption{\it  Contour plot of $\mathfrak{d}$ for $\ell=0.4$ mm, circular inclusion and $f=1/4$. (a)  $\bar{H}_{22}=0.014 $  (b) $\bar{H}_{22}=0.0172 $  }
  \label{Fig::Specimen_l0_04}
\end{figure}

Figures \ref{Fig::Specimen_l0_01} and \ref{Fig::Specimen_l0_04} show two contour plots of the phase field variable $\mathfrak{d}$ for $\ell=0.1$ mm and $\ell=0.4$ mm, respectively, at two consecutive pseudo-time steps of the quasi-static computation in the case of a circular inclusion and a volumetric content $f=1/4$.
In particular Fig. \ref{Fig::Specimen_l0_01}-(a) refers to $\bar{H}_{22}=0.0062$, Fig. \ref{Fig::Specimen_l0_01}-(b) to $\bar{H}_{22}=0.0072$ , Fig. \ref{Fig::Specimen_l0_04}-(a) to $\bar{H}_{22}=0.014 $, and  Fig. \ref{Fig::Specimen_l0_04}-(b) to $\bar{H}_{22}=0.0172$, respectively. As expected, the spread of the damaged zone increases by increasing the internal length $\ell$.

\begin{figure}[h!]
  \centering
  \begin{tabular}{c c  c}
 \hspace{-0.75cm}
  \includegraphics[width=4.8cm]{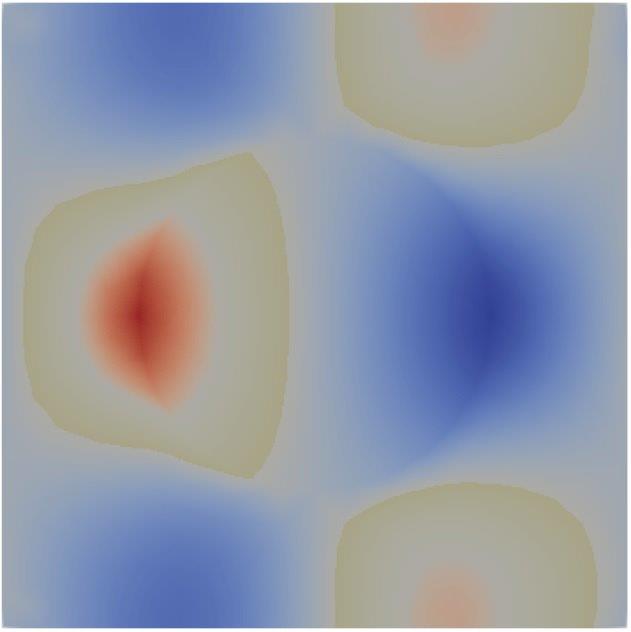}
  &
 \hspace{0.5cm}
  \includegraphics[width=4.8cm]{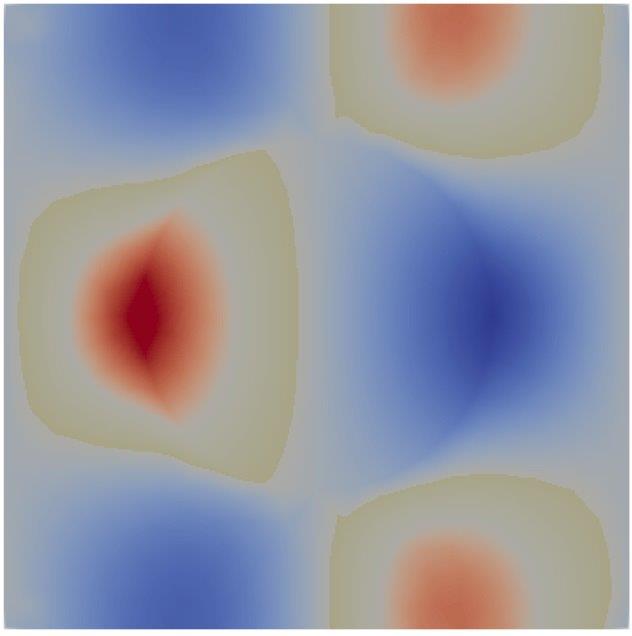}
   &
 \hspace{-0.6cm}
  \includegraphics[width=1.7cm]{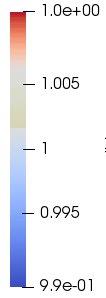}
  \\
   \hspace{-0.75cm} (a)  &  \hspace{0.5cm} (b) &
  \end{tabular}
  \caption{\it Dimensionless micro displacement $\tilde{u}_1(\mathbf{x},\tensor{\xi})$ at point $\mathbf{x}=\{0.75, 0.12\}$ mm of the specimen shown in Fig. \ref{Fig::Specimen}-(b): $\ell=0.4$ mm, circular inclusion with $f=1/4$. (a)  $\bar{H}_{22}=0.014$, (b)  $\bar{H}_{22}=0.0172$.}
  \label{Fig::umicro_ll04}
\end{figure}

\begin{figure}[h!]
  \centering
  \begin{tabular}{c c c }
 \hspace{-0.75cm}
  \includegraphics[width=4.8cm]{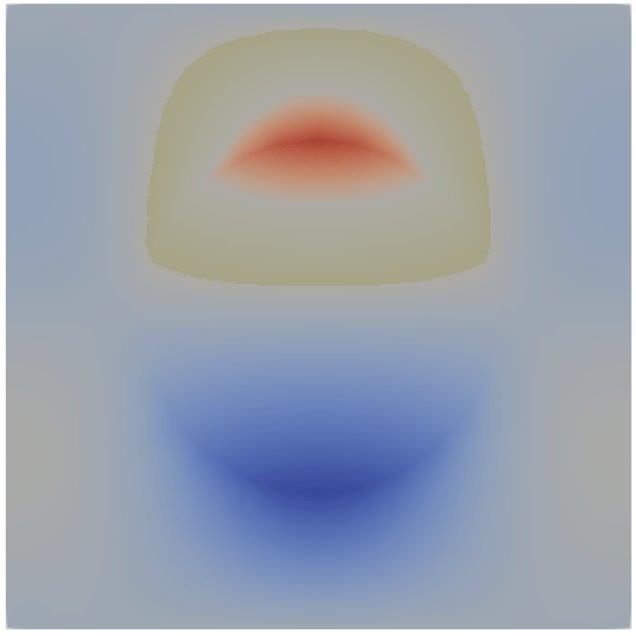}
  &
 \hspace{0.5cm}
  \includegraphics[width=4.8 cm]{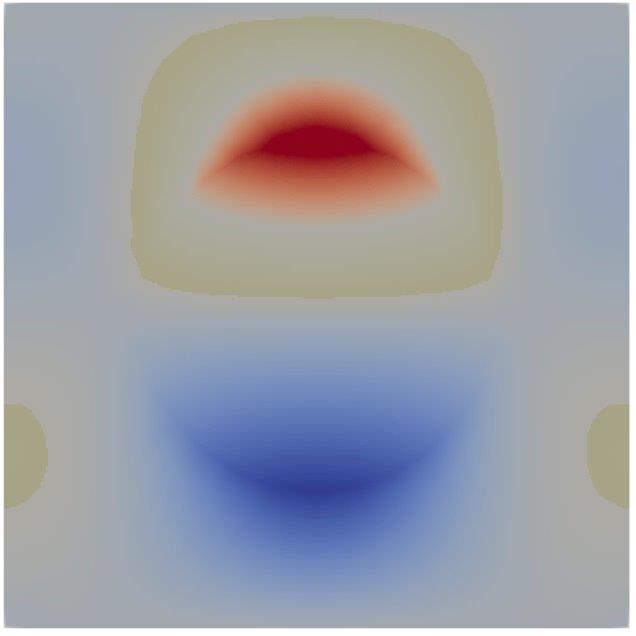}
  &
 \hspace{-0.6cm}
  \includegraphics[width=1.7cm]{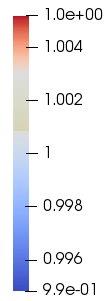}
  \\
   \hspace{-0.75cm} (a)  &  \hspace{0.5cm} (b) &
  \end{tabular}
  \caption{\it Dimensionless micro displacement $\tilde{u}_2(\mathbf{x},\tensor{\xi})$ at point $\mathbf{x}=\{0.75, 0.12\}$ mm of the specimen shown in Fig. \ref{Fig::Specimen}-(b): $\ell=0.4$ mm, circular inclusion and $f=1/4$. (a)  $\bar{H}_{22}=0.014$, (b)  $\bar{H}_{22}=0.0172$.  }
  \label{Fig::vmicro_ll04}
\end{figure}

Without any loss of generality, the cases shown in Fig. \ref{Fig::Specimen_l0_04} have been considered as possible examples to show the applicability of the down-scaling relations established by asymptotic homogenization to reconstruct the local microscopic fields.
In fact, once the macro displacement field $\mathbf{U}(\mathbf{x})$ is known in each point $\mathbf{x}$ of the macroscale model, it is possible to reconstruct the micro displacement field $\mathbf{u}\left(\mathbf{x},\tensor{\xi}\right)$ on the periodic cell by means of the down-scaling relation \eqref{eq:downscaling}, herein truncated to the first order of $\varepsilon$.
For example, considering a point of coordinates $\mathbf{x}=\{0.75, 0.12\}$ mm  inside the specimen sketched in Fig. \ref{Fig::Specimen}-(b), the computed dimensionless micro displacement $\tilde{u}_1(\mathbf{x},\tensor{\xi})=u_1(\mathbf{x},\tensor{\xi})/U_1(\mathbf{x})$ is shown in Fig. \ref{Fig::umicro_ll04} for $\ell=0.4 $ mm and a circular inclusion with $f=1/4$, for two different values of $\bar{H}_{22}$.
In particular, Fig. \ref{Fig::umicro_ll04}-(a) refers to an imposed average strain equal to $\bar{H}_{22}=0.014$, and a macro displacement $U_1(\mathbf{x})=-2.0113\,10^{-3}$ mm, while Fig. \ref{Fig::umicro_ll04}-(b) refers to $\bar{H}_{22}=0.0172$  and $U_1(\mathbf{x})=-2.2932\,10^{-3}$ mm.
The dimensionless micro displacement  $\tilde{u}_2(\mathbf{x},\tensor{\xi})=u_2(\mathbf{x},\tensor{\xi})/U_2(\mathbf{x})$ in the $\mathbf{e}_2$ direction, at the same point $\mathbf{x}=\{0.75, 0.12\}$ mm, is shown in Fig. \ref{Fig::vmicro_ll04}.
Fig. \ref{Fig::vmicro_ll04}-(a) refers to an imposed strain equal to $\bar{H}_{22}=0.014$  and $U_2(\mathbf{x})=-1.1741\,10^{-2}$ mm, while Fig.\ref{Fig::vmicro_ll04}-(b) refers to  $\bar{H}_{22}=0.0172$ and $U_2(\mathbf{x})=-1.4217\,10^{-2}$ mm.

\newpage
\section{Conclusions}
\label{Sec::Conclusions}
The evolution of damage in  composite materials with periodic or quasi-periodic microstructure has been  investigated in the present work trough a  multiscale finite element-based  computational approach, integrating the phase field method with an  homogenization scheme. Specifically, the overall constitutive properties of the heterogeneous material have been determined in closed form by means of a two-scale asymptotic homogenization technique and the  damage  propagation inside  the equivalent homogeneous medium  has been described at the macroscale   via a phase field approach.

The current technique allowed triggering macroscopic damage-like by means of tracking the matrix degradation at the microscale. In this regard, differing from alternative FE$^2$-techniques, this method provided a fast and efficient computational tool for  the  link between the two-scales using the concept of a look-up table scheme. Therefore, it was precluded the instantaneous computation of both scales throughout the simulation.

The performance of the  present technique has been demonstrated by means of several representative  examples. These results clearly evidenced that, for cases under analysis, the overall response of the specimen was   affected by the material length scale, which is usually related to the apparent material strength of the material, and by the micromechanical information regarding the inclusion shape and volume content as well as the damage extent into the matrix. Moreover, through these numerical  applications, it has been shown that these aspects played a significant role in the  the peak response of the specimen and the posterior softening evolution upon failure, featuring a progressive softening instead of an abrupt drop, as in most of the phase field methods of fracture.

Finally, it is worth to recall  that the influence of the microscale over the macroscale results  is one of the  major novelties of the present approach in comparison to  standard phase field methods  for homogeneous materials with no specific scale separation.

\section*{Acknowledgements}
AB would like to acknowledge the financial support by National Group of Mathematical Physics (GNFM-INdAM).
MP would like to acknowledge the financial support of the Italian Ministry of Education, University and Research to the Research Project of National Interest (PRIN 2017).
The authors would like to thank the IMT School for Advanced Studies Lucca for its support to the stays of FF and JR in the IMT Campus as visiting researchers in 2019, making possible the realization of this joint work.


\bigskip
\bibliographystyle{elsarticle-harv}
\bibliography{Bibliography}

\begin{appendices}
\section{Finite element formulation  of the coupled model}
\label{Appendix:finite element}
This appendix details the finite element formulation that has been exploited and implemented in the finite element software FEAP in order to solve the coupled system \eqref{eq:Euler-Lagrange equations} in terms of the macro displacements $U_i(\mathbf{x})$ and the phase field  variable $\mathfrak{d}(\mathbf{x})$.

The weak form of balance equations \eqref{eq:Euler-Lagrange equations} of the coupled field problem detailed in Section \ref{Sec::phase_field}, taking into account boundary conditions \eqref{eq:Neumann_boundary_conditions}, reads
\begin{eqnarray}
\label{eq:weak_form_appendix}
&&- \int_{\Omega}\frac{\partial\psi_{U_h}}{\partial x_k}\,C_{ijhk}\left(\mathfrak{d}\right)\frac{\partial U_i}{\partial x_j}\,\text{d}\Omega
+
\int_{\Gamma_t} \psi_{U_h}\,t_h \,\text{d}\Gamma + \int_{\Omega}\psi_{U_h}\, b_h \, \text{d}\Omega =0
\hspace{0.2cm}\forall
\psi_{U_h}\hspace{0.1cm}s.t.\hspace{0.1cm}\psi_{U_h}=0\hspace{0.1cm}on\hspace{0.1cm}\Gamma_u
\nonumber\\
&&\int_{\Omega}\ell^2\frac{\partial\psi_{\mathfrak{d}}}{\partial x_j}\,\frac{\partial\mathfrak{d}}{\partial x_j}\,\text{d}\Omega +
\int_{\Omega}\psi_{\mathfrak{d}}\,\mathfrak{d}\,\text{d}\Omega +
\int_{\Omega}\frac{\ell}{2 G_C}\psi_\mathfrak{d} H_{ij}\,\frac{\partial C_{ijhk}(\mathfrak{d})}{\partial\mathfrak{d}}H_{hk}\,\text{d}\Omega=0\hspace{0.2cm}
\forall
\psi_{\mathfrak{d}}
\end{eqnarray}
where $\psi_{U_h}$ and $\psi_{\mathfrak{d}}$ are taken as test functions.
Considering the finite dimensional space $V_h$, for which $\{N_j|j=1,2,...,N_h\}$ is a basis, in the finite element discretization the macro displacement field $\mathbf{U}(\mathbf{x})$ and the phase field $\mathfrak{d}(\mathbf{x})$ are approximated as linear combinations of shape functions $N_j(\mathbf{x})$ and nodal unknowns $U_{ij}$ and $\mathfrak{d}_j$
\begin{equation}
\label{eq:discretization}
U_i(\mathbf{x})=\sum_{j=1}^{N_h}N_j(\mathbf{x})U_{ij},
\hspace{0.5cm}
\mathfrak{d}(\mathbf{x})=\sum_{j=1}^{N_h}N_j(\mathbf{x})\mathfrak{d}_{j}
\end{equation}
Analogous approximations are considered for test functions $\psi_{U_i}$ and $\psi_{\mathfrak{d}}$, whose nodal unknowns are indicated respectively as $\delta U_{ij}$ and $\delta\mathfrak{d}_j$
\begin{equation}
\label{eq_dicretization_test_function}
\psi_{U_i} (\mathbf{x})=\sum_{j=1}^{N_h}N_j(\mathbf{x})\delta U_{ij},
\hspace{0.5cm}
\psi_\mathfrak{d}(\mathbf{x})=\sum_{j=1}^{N_h}N_j(\mathbf{x})\delta \mathfrak{d}_{j}
\end{equation}
In every single finite element $e$, one can define matrices $\mathbf{B}_U$ and $\mathbf{B}_{\mathfrak{d}}$ as
\begin{equation}
\label{eq:B_matrices}
\mathbf{B}_U=\mathbf{D}_U\mathbf{N}_U,
\hspace{0.5cm}
\mathbf{B}_\mathfrak{d}=\mathbf{D}_\mathfrak{d}\mathbf{N}_\mathfrak{d}
\end{equation}
where, in a two dimensional setting, matrices $\mathbf{D}_U$ and $\mathbf{D}_{\mathfrak{d}}$ contain the derivatives with respect to coordinates $x_1$ and $x_2$
\begin{equation}
\label{eq:D_matrices}
\mathbf{D}_U
=
\left[
\begin{array}{c c}
\frac{\partial}{\partial x_1} & 0 \\
0 & \frac{\partial}{\partial x_2} \\
\frac{\partial }{\partial x_2} & \frac{\partial }{\partial x_1}
\end{array}
\right],
\hspace{0.5cm}
\mathbf{D}_{\mathfrak{d}}=
\left[
\begin{array}{c}
\frac{\partial}{\partial x_1}\\
\frac{\partial}{\partial x_2}
\end{array}
\right],
\end{equation}
while $\mathbf{N}_U$ and $\mathbf{N}_{\mathfrak{d}}$ collect the shape functions
\begin{equation}
\label{eq:N_matrices}
\mathbf{N}_U
=
\left[
\begin{array}{c c c c c c c}
N_1 & 0 & N_2 & 0 & ... & N_{Nnod} & 0\\
0 & N_1 & 0 & N_2 & ... & 0 & N_{Nnod}
\end{array}
\right],
\hspace{0.5cm}
\mathbf{N}_{\mathfrak{d}}=
\left[
\begin{array}{c c c c}
N_1 & N_2 & ... & N_{Nnod}
\end{array}
\right],
\end{equation}
being $N_{nod}$ the number of single element nodes.
Thus, the weak form \eqref{eq:weak_form_appendix} can be written over each element domain $\Omega_e$ as
\begin{eqnarray}
\label{eq:elemental_matrix_form}
&&-  \tensor{\delta}\mathbf{U}^T
\int_{\Omega_e}
\mathbf{B}_U^T\,
\mathbf{C}(\mathfrak{d})
\,
\mathbf{B}_U
\,\text{d}\Omega\,\mathbf{U}+\tensor{\delta}\mathbf{U}^T
\int_{\Gamma_{e_t}}
\mathbf{N}_U^T\,
\mathbf{t}
\,\text{d}\Gamma\,
+
\tensor{\delta}\mathbf{U}^T
\int_{\Omega_e}
\mathbf{N}_U^T\,
\mathbf{b}
\,\text{d}\Omega=0,
\nonumber\\
&&
\tensor{\delta}\mathfrak{d}^T
G_C \ell
\int_{\Omega_e}\mathbf{B}_{\mathfrak{d}}^T\, \mathbf{B}_{\mathfrak{d}}\,\text{d}\Omega\,\mathfrak{d}
+
\tensor{\delta}\mathfrak{d}^T
\frac{G_C}{\ell}\int_{\Omega_e}\mathbf{N}_\mathfrak{d}^T
\mathbf{N}_{\mathfrak{d}}\,\text{d}\Omega\,\mathfrak{d}
+
\frac{1}{2}\mathbf{\delta}\mathfrak{d}^T
\int_{\Omega_e}\mathbf{N}_{\mathfrak{d}}^T \mathbf{U}^T\mathbf{B}_U^T\,\frac{\partial \mathbf{C}(\mathfrak{d})}{\partial\mathfrak{d}}\mathbf{B}_U\,\text{d}\Omega\, \mathbf{U}
= 0
\end{eqnarray}
The numerical solution of the coupled problem \eqref{eq:elemental_matrix_form} is obtained by means of an iterative Newton-Raphson procedure, for which residual vectors of the displacement field $\mathbf{R_U}$ and of the phase field $\mathbf{R}_{\mathfrak{d}}$ are defined as
\begin{subequations}
\begin{align}
&\mathbf{R}_U=-\int_{\Omega_e}
\mathbf{B}_U^T\,
\mathbf{C}(\mathfrak{d})
\,
\mathbf{B}_U
\,\text{d}\Omega\,\mathbf{U}
 +
\int_{\Gamma_{e_t}}
\mathbf{N}_U^T\,
\mathbf{t}
\,\text{d}\Gamma\,
+
\int_{\Omega_e}
\mathbf{N}_U^T\,
\mathbf{b}
\,\text{d}\Omega,\\
&\mathbf{R}_\mathfrak{d}
=
-G_C\ell
\int_{\Omega_e}\mathbf{B}_{\mathfrak{d}}^T   \mathbf{B}_{\mathfrak{d}}\,\text{d}\Omega\,\mathfrak{d}
-
\frac{G_C}{\ell}
\int_{\Omega_e}\mathbf{N}_\mathfrak{d}^T
\mathbf{N}_{\mathfrak{d}}\,\text{d}\Omega\,\mathfrak{d}
-
\frac{1}{2}
\int_{\Omega_e}\mathbf{N}_{\mathfrak{d}}^T \mathbf{U}^T\mathbf{B}_U^T\,\frac{\partial \mathbf{C}(\mathfrak{d})}{\partial\mathfrak{d}}\mathbf{B}_U\,\text{d}\Omega\, \mathbf{U}
\end{align}
\label{eq:residuals}
\end{subequations}
The specific form of elemental stiffness matrices reads
\begin{subequations}
\label{eq:stiffness_matrices}
\begin{align}
&\mathbf{K}_{UU}^e=-\frac{\partial \mathbf{R}_U}{\partial\mathbf{U}}=\int_{\Omega_e}
\mathbf{B}_U^T\,
\mathbf{C}(\mathfrak{d})
\,
\mathbf{B}_U
\,\text{d}\Omega,\\
&\mathbf{K}_{U\mathfrak{d}}^e=-\frac{\partial\mathbf{R}_U}{\partial \mathfrak{d}}=
\int_{\Omega_e}
\mathbf{B}_U^T\,
\frac{\partial \mathbf{C}(\mathfrak{d})}{\partial \mathfrak{d}}
\,
\mathbf{B}_U\,\mathbf{U}\,\mathbf{N}_{\mathfrak{d}}
\,\text{d}\Omega,\\
&\mathbf{K}_{\mathfrak{d}U}^e=-\frac{\partial\mathbf{R}_{\mathfrak{d}}}{\partial U}=
\int_{\Omega_e}
\mathbf{N}_{\mathfrak{d}}^T
\,\mathbf{U}^T\,\mathbf{B}_U^T\,\frac{\partial\mathbf{C}(\partial\mathfrak{d})}{\partial\mathfrak{d}}\mathbf{B}_U\,\text{d}\Omega,\\
&\mathbf{K}_{\mathfrak{d}\mathfrak{d}}^e=-
\frac{\partial \mathbf{R}_{\mathfrak{d}}}{\partial\mathfrak{d}}=G_C\ell\int_{\Omega_e}
\mathbf{B}_{\mathfrak{d}}^T\,\mathbf{B}_{\mathfrak{d}}\,\text{d}\Omega
+\frac{G_C}{\ell}\,
\int_{\Omega_e}
\mathbf{N}_{\mathfrak{d}}^T\,
\mathbf{N}_{\mathfrak{d}}\,\text{d}\Omega
+
\frac{1}{2}
\int_{\Omega_e}
\mathbf{N}_{\mathfrak{d}}^T
\,
\mathbf{U}^T\mathbf{B}_U^T
\frac{\partial^2 \mathbf{C}(\mathfrak{d})}{\partial \mathfrak{d}^2}
\mathbf{B}_U
\mathbf{U}\,\mathbf{N_{\mathfrak{d}}}\,\text{d}\Omega
\end{align}
\end{subequations}
Consistently with the linearization of the resulting nonlinear system of equations \eqref{eq:elemental_matrix_form}, at each iteration of the Newton-Raphson loop, the following linear system has to be solved
\begin{equation}
\label{eq:linear_system_NR}
\left[
\begin{array}{c c}
\mathbf{K}_{UU} & \mathbf{K}_{U\mathfrak{d}} \\
\mathbf{K}_{\mathfrak{d} U} & \mathbf{K}_{\mathfrak{d}\mathfrak{d}}
\end{array}
\right]
\left[
\begin{array}{c}
\Delta \mathbf{U}\\
\Delta \mathfrak{d}
\end{array}
\right]
=
\left[
\begin{array}{c}
\mathbf{R}_U\\
\mathbf{R}_{\mathfrak{d}}
\end{array}
\right]
\end{equation}
where $\Delta\mathbf{U}$ and $\Delta\mathfrak{d}$ are discretized according to \eqref{eq:discretization} and the elemental stiffness matrices \eqref{eq:stiffness_matrices} and the residual vectors \eqref{eq:residuals} have been assembled in the corresponding global ones\footnote{Finally note that following \citep{miehe2010}, the current formulation is equipped with a viscous crack resistance parameter, leading to the modification of the operators associated with the phase field variable. An alternative solution scheme would encompass a staggered Jacobi-type method which can be easily recalled by eliminating the coupling stiffness matrices and adopting an alternate minimization procedure}.

Coupling with asymptotic homogenization is established by taking the closed-form expression for $\mathbf{C}(\mathfrak{d})$ and $\partial\mathbf{C}(\mathfrak{d})/\partial\mathfrak{d}$ provided by an off-line computation based on  asymptotic homogenization, for different values of $\mathfrak{d}$.

\end{appendices}
\end{document}